\documentclass{article}

\usepackage{amssymb}
\usepackage{amsmath}
\usepackage{amsfonts}
\usepackage{epsfig}
\usepackage[usenames]{color}

\newcommand{\Au}{{\bf (A.1)}}
\newcommand{\Ad}{{\bf (A.2)}}
\newcommand{\At}{{\bf (A.3)}}
\newcommand{\Aq}{{\bf (A.4)}}

\newcommand{\muu}{\mu_x^{(1)}(h)}
\newcommand{\mud}{\mu_x^{(2)}(h)}

\newcommand{\baF}{\bar{F}}
\newcommand{\habaF}{\hat{\bar{F_n}}}
\newcommand{\invhabaF}{\hat{\bar{F_n}}^{\leftarrow}}
\newcommand{\toP}{\stackrel{P}{\longrightarrow}}
\newcommand{\tod}{\stackrel{d}{\longrightarrow}}

\newcommand{\R}{\mathbb{R}}
\newcommand{\PP}{\mathbb{P}}
\newcommand{\I}{\mathbb{I}}
\newcommand{\E}{\mathbb{E}}

\newcommand{\gampick}{\hat\gamma_n^{\phi_{\mbox{\tiny FP}}}(x)}
\newcommand{\gamhill}{\hat\gamma_n^{\phi_1}(x)}
\newcommand{\qweis}{\hat q_n^{\mbox{\tiny W}}(\beta_n|x)}

\newcommand{\CQFD}
{%
\mbox{}%
\nolinebreak%
\hfill%
\rule{2mm}{2mm}%
\medbreak%
\par%
}
 \newcommand{\proof}{{\bf Proof. }}

\newtheorem{Theo}{Theorem}

\newtheorem{Coro}{Corollary}
\newtheorem{Lem}{Lemma}

\begin{document}

\title{Functional kernel estimators of large conditional quantiles}
\author{Laurent Gardes$^{(1)}$ \& St\'ephane Girard$^{(2,\star)}$  }
\date{$^{(1)}$
IRMA, Universit\'e de Strasbourg,
7 rue Ren\'e Descartes,\\
67084 Strasbourg Cedex, France.\\
$^{(2)}$
INRIA Rh\^one-Alpes \& LJK, team Mistis, Inovall\'ee, 655, av. de l'Europe, Montbonnot, 38334 Saint-Ismier cedex, France.
\\
$^{(\star)}$  {\tt Stephane.Girard@inria.fr} (corresponding author)
}
\maketitle
\begin{abstract}
We address the estimation of conditional quantiles 
when the covariate is functional and when the order of 
the quantiles converges to one as the sample size increases.
In a first time, we investigate to what extent these large conditional 
quantiles can still be estimated through a functional kernel estimator
of the conditional survival function. 
Sufficient conditions on the rate of convergence of their order to one
are provided to obtain asymptotically Gaussian distributed estimators.
In a second time, basing on these result, a functional Weissman estimator is derived, permitting
to estimate large conditional quantiles of arbitrary large order.
These results are illustrated on finite sample situations.\\


\end{abstract}

\section{Introduction}
\label{intro}
Let $(X_i,Y_i)$, $i=1,\dots,n$ be independent
copies of a random pair $(X,Y)$ in $E \times\R$
where $E$ is an infinite dimensional space associated to a semi-metric $d$.
We address the problem of estimating 
$q(\alpha_n|x)\in\R$
verifying 
$
\PP(Y>q(\alpha_n|x)|X=x)=\alpha_n
$
where $\alpha_n\to 0$ as $n\to\infty$ and $x\in E$.
In such a case, $q(\alpha_n|x)$ is referred to as a large
conditional quantile in contrast to classical conditional quantiles
(or regression quantiles)
for which $\alpha_n=\alpha$ is fixed in $(0,1)$.
While the nonparametric estimation of ordinary regression quantiles 
has been extensively studied 
(see for instance~\cite{Roussas,Stone} or~\cite{LivreToulouse}, Chapter~5),
less attention has been paid to large conditional quantiles despite their potential
interest. 
In climatology, large conditional quantiles may explain how climate change over years might affect extreme temperatures. 
In the financial econometrics literature, they illustrate the link
between extreme hedge fund returns and some measures of risk. 
Parametric models are introduced
in~\cite{davsmi90,smith89} and semi-parametric methods are
considered in~\cite{beigoe03,haltaj00}.
Fully non-parametric estimators have been first introduced in~\cite{davram00, chadav05}
through local polynomial and spline models.
In both cases, the authors focus on univariate covariates and on the
finite sample properties of the estimators.
Nonparametric methods based on moving windows and nearest neighbors are introduced 
respectively in~\cite{GG} and~\cite{Extremes2}.
We also refer to~\cite{FHR}, Theorem~3.5.2,
for the approximation of the
nearest neighbors distribution using the Hellinger distance and to~\cite{Gango}
for the study of their asymptotic distribution.

An important literature is devoted to the particular case where 
the conditional distribution of $Y$ given $X=x$ has a finite endpoint 
$\varphi(x)$ and when $X$ is a finite dimensional random variable.
The function $\varphi$ is referred to as the frontier
and can be estimated from an estimator of the conditional quantile
$q(\alpha_n|x)$ with $\alpha_n\to 0$.
As an example, a kernel estimator of $\varphi$ is proposed in~\cite{JMVA2},
the asymptotic normality being proved only
when $Y$ given $X=x$ is uniformly distributed
on $[0,\varphi(x)]$. 
We refer to~\cite{KorTsy} for a review on this topic.

Estimation of unconditional large quantiles is also widely studied
since the introduction of Weissman estimator~\cite{wei78} dedicated
to heavy-tailed distributions, Weibull-tail estimators~\cite{WT2,WT1} dedicated to
light-tailed distributions and Dekkers and de Haan estimator~\cite{Dek}
adapted to the general case.

In this paper, we focus on the setting where the
conditional distribution of $Y$ given $X=x$ has an infinite endpoint
and is heavy-tailed,
an analytical characterization of this property being given in the next section.
In such a case, the frontier function does not exist and
$q(\alpha_n|x)\to\infty$ as $\alpha_n\to 0$.
Nevertheless, we show, under some conditions, that large regression
quantiles $q(\alpha_n|x)$ can still be estimated through a functional kernel estimator
of $\PP(Y>.|x)$. We provide sufficient conditions on the rate
of convergence of $\alpha_n$ to 0 so that our estimator is
asymptotically Gaussian distributed.
Making use of this, some functional estimators of the conditional tail-index
are introduced and a functional Weissman estimator~\cite{wei78} is derived, permitting
to estimate large conditional quantiles $q(\beta_n|x)$ where $\beta_n\to 0$
arbitrarily fast.

Assumptions are introduced and discussed in Section~\ref{hypo}.
Main results are provided in Section~\ref{main} and
illustrated both on simulated and real data in Section~\ref{simul} and Section~\ref{realdata}.
Extensions of this work are briefly discussed in Section~\ref{further}.
Proofs are postponed to the appendix.

\section{Notations and assumptions}
\label{hypo}

The conditional survival function (csf) of $Y$ given $X=x$ is
denoted by $\baF(y|x)=\PP(Y>y|X=x)$. 
The functional estimator of $\baF(y|x)$ is defined for all $(x,y)\in E \times\R$
by
\begin{equation}
\label{defestproba}
\habaF(y|x)= \left.\sum_{i=1}^n K(d(x,X_i)/h) Q((Y_i-y)/\lambda) \right/ 
\sum_{i=1}^n K(d(x,X_i)/h) ,
\end{equation}
with $Q(t)=\int_{-\infty}^t Q'(s) ds$ where 
$K:\R^+\to\R^+$ and $Q':\R\to\R^+$ are two kernel functions, and $h=h_n$ and $\lambda=\lambda_n$ are two nonrandom sequences
(called window-width) such that 
$h\to 0$ as $n\to\infty$. Let us emphasize that the condition $\lambda\to 0$ is not required in this context.
This estimator was considered for instance in~\cite{LivreToulouse}, page~56.
Its rate of uniform strong consistency is established by~\cite{Ferraty2}.
In Theorem~\ref{thproba} hereafter, the asymptotic distribution of~(\ref{defestproba})
is established when estimating small tail probabilities, {\it i.e}
when $y=y_n$ goes to infinity with the sample size $n$.
Similarly, the functional estimators of conditional quantiles 
$q(\alpha|x)$
are defined via the generalized inverse of $\habaF(.|x)$: 
\begin{equation}
\label{defestquant}
\hat q_n(\alpha | x ) = \invhabaF(\alpha|x)=\inf\{t,\; \habaF(t|x)\leq \alpha\},
\end{equation}
for all $\alpha\in(0,1)$. 
Many authors are interested in this estimator for fixed $\alpha\in(0,1)$. 
Weak and strong consistency are proved respectively in~\cite{Stone} and~\cite{Gannoun}.
Asymptotic normality is shown in~\cite{Berlinet, Samanta, Stute} when $E$ is finite dimensional
and by~\cite{Ferraty} for a general metric space under dependence assumptions.
In Theorem~\ref{thquant},
the asymptotic distribution of~(\ref{defestquant})
is investigated when estimating large quantiles, {\it i.e}
when $\alpha=\alpha_n$ goes to 0 as the sample size $n$ goes to infinity.
The asymptotic behavior of such estimators depends on the nature of
the conditional distribution tail. In this paper, we focus on heavy tails.
More specifically, we assume that the csf satisfies\\

\Au: $\displaystyle\baF(y|x)=c(x)  \exp\left\{-\int_1^y \left( \frac{1}{\gamma(x)}-\varepsilon(u|x)\right)\frac{du}{u}\right\}$,\\

\noindent where $\gamma$ is a positive function of the covariate $x$,
$c$ is a positive function and $|\varepsilon(.|x)|$ is continuous and ultimately 
decreasing to 0. Examples of such distributions are provided in Table~\ref{exemples}.
\Au~implies that the conditional distribution
of $Y$ given $X=x$ is in the Fr\'echet maximum domain of attraction.
In this context, $\gamma(x)$ is referred to as the conditional tail-index 
since it tunes the tail heaviness of the conditional distribution of
$Y$ given $X=x$. More details on extreme-value theory can be found for
instance in~\cite{EMBR}.
Assumption \Au~also yields that
$\baF(.|x)$ is regularly varying at infinity with index $-1/\gamma(x)$.
 {\it i.e} for all $\zeta>0$,
\begin{equation}
\label{defl}
 \lim_{y\to\infty} \frac{\baF(\zeta y | x)}{\baF(y|x)} = \zeta^{-1/\gamma(x)}.
\end{equation}
We refer to~\cite{BING} for a general account on regular variation theory.
The auxiliary function $\varepsilon(.|x)$ plays
an important role in extreme-value theory since it drives the
speed of convergence in~(\ref{defl}) and more generally the bias
of extreme-value estimators. Therefore, it may be of interest to
specify how it converges to~0. In~\cite{alves2, Gomes1},
$|\varepsilon(.|x)|$ is supposed to be regularly varying
and the estimation of the corresponding regular variation index is addressed.

\noindent Some Lipschitz conditions are also required:

\Ad: There exist $\kappa_\varepsilon$, $\kappa_c$, $\kappa_\gamma>0$ and $u_0>1$ such that for all $(x,x')\in E\times E$ and $u>u_0$,
\begin{eqnarray*}
\displaystyle \left|\log c(x)-\log c(x')\right|
&\leq& \kappa_c d(x,x'),\\
\displaystyle \left|\varepsilon(u|x)-\varepsilon(u|x')\right|
&\leq& \kappa_\varepsilon d(x,x'),\\
\displaystyle \left|\frac{1}{\gamma(x)}-\frac{1}{\gamma(x')}\right|
&\leq &\kappa_\gamma d(x,x').
\end{eqnarray*}

\noindent The last two assumptions are standard in the functional kernel estimation framework.

\At: $K$ is a function with support $[0,1]$
and there exist $0<C_1<C_2<\infty$ such that $C_1\leq K(t) \leq C_2$ for all $t\in[0,1]$.

\Aq: $Q'$ is a probability density function (pdf) with support $[-1,1]$.

\noindent One may also assume without loss of generality that $K$ integrates to one.
In this case, 
$K$ is called a type I kernel, see~\cite{LivreToulouse}, Definition~4.1.
Letting $B(x,h)$ be the ball of center $x$ and radius $h$,
we finally introduce $\varphi_x(h):=\PP(X\in B(x,h)$ the small ball probability of $X$.
Under \At, the $\tau$-th moment $\mu^{(\tau)}_x(h):=\E\{K^\tau(d(x,X)/h)$ can be controlled for all $\tau>0$ by Lemma~\ref{lemmu} in 
Appendix. It is shown that $\mu^{(\tau)}_x(h)$ is of the same asymptotic order as $\varphi_x(h)$.

\section{Main results}
\label{main}

The first step towards the estimation of large conditional quantiles
is the estimation of small tail probabilities
$\baF(y_n|x)$ when $y_n\to\infty$ as $n\to\infty$. 

\subsection{Estimation of small tail probabilities}

Defining 
$$
\Lambda_n(x)= \left(n  \baF(y_n|x) \frac{(\muu)^2}{\mud} \right)^{-1/2},
$$
the following result provides sufficient conditions for 
the asymptotic normality of $\habaF(y_n|x)$.

\begin{Theo}
\label{thproba}
Suppose \Au~-- \Aq~hold. Let $x\in E$ such that $\varphi_x(h)>0$ and  introduce $y_{n,j}=a_j y_n(1+o(1))$ for $j=1,\dots, J$ with
$0<a_1<a_2< \dots < a_J$ and where $J$ is a positive integer.
If $y_n\to\infty$ such that  $n\varphi_x(h)  \baF(y_n|x) \to \infty$ and
$n\varphi_x(h)  \baF(y_n|x) (\lambda/y_n \vee h \log y_n)^2  \to 0$ as $n\to \infty$, then
$$
\left\{ \Lambda_n^{-1}(x)\left(\frac{\habaF(y_{n,j}|x)
}{\baF(y_{n,j}|x)}-1\right)\right\}_{j=1,\dots,J}
$$
is asymptotically Gaussian, centered, with covariance matrix 
$C(x)$ where $C_{j,j'}(x) = a^{1/\gamma(x)}_{j\wedge j'}$ for $(j,j')\in \{1,\dots, J\}^2$.
\end{Theo}
Note that $n \varphi_x(h) \baF(y_n|x) \to \infty$
is a necessary and sufficient condition for 
the almost sure presence of at least one sample point
in the region $B(x,h)\times (y_n,\infty)$ of $E \times \R$,
see Lemma~\ref{leminter} in Appendix. Thus, this natural condition
states that one cannot estimate small tail probabilities out of the sample
using $\habaF$. Besides, from Lemma~\ref{lemmu}, $\Lambda_n^{-2}(x)$ is of the same asymptotic order as
$n \varphi_x(h) \baF(y_n|x)$ and consequently $\Lambda_n(x)\to 0$ as $n\to\infty$.
Theorem~\ref{thproba} thus entails
$
\habaF(y_{n,j}|x)/\baF(y_{n,j}|x)\toP 1
$
which can be read as a consistency of the estimator.
The second condition 
 $n\varphi_x(h)  \baF(y_n|x) (\lambda/y_n \vee h \log y_n)^2  \to 0$
imposes to the biases $\lambda/y_n$ and $h \log y_n$  introduced by the two smoothings to be negligible compared 
to the standard deviation $\Lambda_n(x)$ of the estimator.
Theorem~\ref{thproba} may be compared to~\cite{einmahl} which
establishes the asymptotic behavior of the empirical survival function
in the unconditional case but without assumption on the distribution.

\subsection{Estimation of large conditional quantiles within the sample}

In this paragraph, we focus on the estimation of large conditional quantiles
of order $\alpha_n$ such that $n \varphi_x(h) \alpha_n\to\infty$ as $n\to\infty$.
This is a necessary and sufficient condition for 
the almost sure presence of at least one sample point
in the region $B(x,h)\times (q(\alpha_n|x),\infty)$ of $E \times \R$,
see Lemma~\ref{leminter} in Appendix.
In other words, the large conditional quantile $q(\alpha_n|x)$ is located within the sample.
Letting
$$
\sigma_n(x)=\left(n  \alpha_n \frac{(\muu)^2}{\mud} \right)^{-1/2},
$$
Lemma~\ref{lemmu} shows that $\sigma_n(x)$ is of the same asymptotic order as
$(n \varphi_x(h) \alpha_n)^{-1/2}$ and thus the condition $n \varphi_x(h) \alpha_n\to\infty$
is equivalent to $\sigma_n(x)\to 0$ as $n\to\infty$.
\begin{Theo}
\label{thquant}
Suppose \Au~-- \Aq~hold. Let $x\in E$ such that $\varphi_x(h)>0$ and consider a sequence
$\tau_1>\tau_2> \dots > \tau_J>0$ where $J$ is a positive integer.
If $\alpha_n\to 0$ such that $\sigma_n(x) \to 0$ and $\sigma_n^{-1}(x)( \lambda/q(\alpha_n|x) \vee h \log\alpha_n)\to 0$ 
as $n\to \infty$,
then
$$
\left\{ \sigma_n^{-1}(x)\left(\frac{\hat q_n(\tau_j \alpha_{n}|x)
}{q(\tau_j \alpha_{n}|x)}-1\right)\right\}_{j=1,\dots,J}
$$
is asymptotically Gaussian, centered, with covariance matrix 
$
{\gamma^2(x)} \Sigma
$
where $\Sigma_{j,j'}= 1/\tau_{j\wedge j'}$ for $(j,j')\in \{1,\dots, J\}^2$.
\end{Theo}
Remark that \Au~provides an asymptotic expansion of the density function of $Y$ given $X=x$:
$$
f(y|x)=\frac{1}{\gamma(x)}\frac{\baF(y|x)}{y}(1-\varepsilon(y|x)) = \frac{1}{\gamma(x)}\frac{\baF(y|x)}{y} (1+o(1))
$$
as $y\to\infty$.
Consequently, Theorem~\ref{thquant} entails that the random vector
$$
\left\{ \frac{\muu}{( \mud)^{1/2}} (n \tau_j \alpha_n(1-\tau_j\alpha_n))^{-1/2}f(q(\tau_j \alpha_{n}|x)|x) \left(
\hat q_n(\tau_j \alpha_{n}|x)-q(\tau_j \alpha_{n}|x)\right)\right\}_{j=1,\dots,J}
$$
is also asymptotically Gaussian and centered. This result coincides with~\cite{Berlinet}, Theorem~6.4
established in the case where $\alpha_n=\alpha$ is fixed in $(0,1)$ and in a finite dimensional setting. 

\subsection{Estimation of arbitrary large conditional quantiles}

This paragraph is dedicated to the estimation of large conditional quantiles
of arbitrary small order $\beta_n$.
For instance, if $n \varphi_x(h) \beta_n\to c \in [1,\infty)$ then $q(\beta_n|x)$
is located near the boundary of the sample. If 
$n \varphi_x(h) \beta_n\to c \in [0,1)$ then $q(\beta_n|x)$
is located outside the sample. 
Here, a functional Weissman
estimator~\cite{wei78} is proposed to tackle all possible situations:
\begin{equation}
\label{defweis}
 \qweis= \hat q_n(\alpha_n|x) (\alpha_n/\beta_n)^{\hat\gamma_n(x)}.
\end{equation}
Here, $\hat q_n(\alpha_n|x)$ is the functional estimator~(\ref{defestquant}) of a
large conditional quantile $q(\alpha_n|x)$ within the sample 
and $\hat\gamma_n(x)$ is an estimator of the conditional
tail-index $\gamma(x)$. As illustrated in the next theorem, the extrapolation
factor $(\alpha_n/\beta_n)^{\hat\gamma_n(x)}$ allows to estimate arbitrary large quantiles.
\begin{Theo}
\label{thweis}
Suppose \Au~-- \Aq~hold.
Let $x\in E$ and introduce 
\begin{itemize}
\item $\alpha_n\to 0$ such that $\sigma_n(x) \to 0$ and $ \sigma_n^{-1}(x) (\lambda/q(\alpha_n|x) \vee h\log\alpha_n \vee \varepsilon(q(\alpha_n|x)|x)) \to 0$ as $n\to \infty$,
\item $(\beta_n)$ such that $\beta_n/\alpha_n\to 0$ as $n\to \infty$,
\item $\hat\gamma_n(x)$ such that $\sigma_n^{-1}(x)(\hat\gamma_n(x)-\gamma(x))\tod {\mathcal {N}}(0,V(x))$ where $V(x)>0$.
\end{itemize}
Then, 
$$
\frac{\sigma_n^{-1}(x)}{\log(\alpha_n/\beta_n)}\left(\frac{\qweis} {q(\beta_n|x)}-1\right)\tod {\mathcal{N}}(0,V(x)).
$$
\end{Theo}
\noindent Let us now focus on the estimation of the conditional tail-index.
Let $\alpha_n\to 0$ and consider a sequence $1=\tau_1>\tau_2> \dots > \tau_J>0$ where $J$ is a positive integer.
Two additional notations are introduced for the sake of simplicity:
$u=(1,\dots,1)^t \in \R^J$ and $v=(\log(1/\tau_1),\dots,\log(1/\tau_J))^t \in \R^J$.
The following family of estimators is proposed
\begin{equation}
\label{defestgamma}
\hat\gamma^{\phi}_n(x) = \frac{\phi (\log \hat q_n(\tau_1\alpha_n|x),\dots,\log \hat q_n(\tau_J\alpha_n|x))}
{\phi(\log(1/\tau_1),\dots,\log(1/\tau_J))},
\end{equation}
where $\phi:\R^J\to\R$ denotes a twice differentiable function
verifying the shift and location invariance conditions 
\begin{equation}
\label{condphi}
\left\{
\begin{array}{lll}
\phi( \theta v)&=&\theta \phi(v)\\
\phi( \eta u + x)&=& \phi(x)
\end{array}
\right.
\end{equation}
for all $\theta> 0$, $\eta\in\R$ and $x\in\R^J$. In the case where $J=3$, $\tau_1=1$, $\tau_2=1/2$ and $\tau_3=1/4$, the function
$$
\phi_{\mbox{\tiny FP}}(x_1,x_2,x_3)= \log\left(\frac{\exp (4x_2) - \exp (4x_1)}{\exp (4x_3) - \exp (4x_2)}\right)
$$
leads us to a functional version of Pickands estimator~\cite{Pick}:
$$
\gampick=\frac{1}{\log 2}\log\left(
\frac{\hat q_n(\alpha_n|x) - \hat q_n(2\alpha_n|x)}
{\hat q_n(2\alpha_n|x) - \hat q_n(4\alpha_n|x)}\right).
$$
We refer to~\cite{Gijbels} for a different variant
of Pickands estimator in the context where the distribution of $Y$ given $X=x$ has a finite
endpoint.
Besides, introducing the function $m_p(x_1,\dots,x_J)=\sum_{j=1}^J (x_j-x_1)^p$ for all $p>0$ and considering 
$\phi_p(x)=m_p^{1/p}(x)$ gives rise to a functional version of the estimator considered for instance in~\cite{Segers}, example~(a):
$$
\hat\gamma^{\phi_p}_n(x) =  \left(\sum_{j=1}^J \left[\log \hat q_n(\tau_j \alpha_n|x) -\log \hat q_n(\alpha_n|x)\right]^p \left / \sum_{j=1}^J [\log (1/\tau_j)]^p\right.\right)^{1/p}.
$$
As a particular case $\phi_{1}(x)= m_1(x)$
corresponds to a functional version of the Hill estimator~\cite{Hill}:
$$
\gamhill = \sum_{j=1}^J \left[\log \hat q_n(\tau_j \alpha_n|x) -\log \hat q_n(\alpha_n|x)\right] \left / \sum_{j=1}^J \log (1/\tau_j)\right..
$$
More interestingly, if $\{\phi^{(1)},\dots,\phi^{(H)}\}$ is a set of $H$ functions satisfying~(\ref{condphi})
and if $A:\R^H\to\R$ is a homogeneous function of degree~1, then the aggregated function $A(\phi^{(1)},\dots,\phi^{(H)})$
also satisfies~(\ref{condphi}). 
Generalizations of the functional Hill estimator can then be obtained using $H=2$, $A_p(x,y)=x^p y^{1-p}$ and
defining $\phi_{p,q,r}=A_p(\phi_q,\phi_r)=m_q^{p/q} m_r^{(1-p)/r}$:
$$
\hat\gamma^{\phi_{p,q,r}}_n(x) =  \frac{\left(\sum_{j=1}^J \left[\log \hat q_n(\tau_j \alpha_n|x) -\log \hat q_n(\alpha_n|x)\right]^p\right)^{p/q}
\left(\sum_{j=1}^J [\log (1/\tau_j)]^r\right)^{(p-1)/r}}
{\left(\sum_{j=1}^J \left[\log \hat q_n(\tau_j \alpha_n|x) -\log \hat q_n(\alpha_n|x)\right]^r\right)^{(p-1)/r}
\left(\sum_{j=1}^J [\log (1/\tau_j)]^p\right)^{p/q}}.
$$
For instance, the estimator introduced by~\cite{Gomes2}, equation~(2.2) corresponds to the particular function
$\phi_{p,p,1}$ and the estimator of~\cite{Caeiro} corresponds to $\phi_{p,p\theta,p-1}$.

\noindent For an arbitrary function $\phi$, the asymptotic normality of $\hat\gamma^{\phi}_n(x)$ is a consequence of Theorem~\ref{thquant}.
The following result permits to establish the asymptotic normality of the above mentioned estimators in an unified way.
\begin{Theo}
\label{thgamma}
Under the assumptions of~Theorem~\ref{thquant} and if, moreover, 
$\sigma_n^{-1}(x)\varepsilon(q(\alpha_n|x)|x)\to 0$ as $n\to \infty$, then,
$\sigma_n^{-1}(x) (\hat\gamma^{\phi}_n(x)-\gamma(x))$ converges to a centered Gaussian
random variable with variance
$$
V_\phi(x)=\frac{\gamma^2(x)}{\phi^2(v)}(\nabla\phi(\gamma(x)v))^t \Sigma  (\nabla\phi(\gamma(x)v)).
$$
\end{Theo}
Let us note that the additional condition $\sigma_n^{-1}(x)\varepsilon(q(\alpha_n|x)|x)\to 0$
is standard in the extreme-value framework: Neglecting the unknown function $\varepsilon(.|x)$ in
the construction of $\hat\gamma^{\phi}_n(x)$ yields a bias that should be negligible with respect
to the standard deviation $\sigma_n(x)$ of the estimator.
Finally, combining Theorem~\ref{thweis} and Theorem~\ref{thgamma}, the asymptotic distribution
of the functional large quantile estimator $\hat q_n^{\mbox{\tiny W},\phi}(\beta_n|x)$ based on~(\ref{defweis}) and~(\ref{defestgamma}) is readily obtained.
\begin{Coro}
Suppose \Au~-- \Aq~hold. Let $x\in E$ such that $\varphi_x(h)>0$ and consider a sequence
$1=\tau_1>\tau_2> \dots > \tau_J>0$ where $J$ is a positive integer. If
\begin{itemize}
\item $\alpha_n\to 0$, $\sigma_n(x) \to 0$ and $ \sigma_n^{-1}(x) (\lambda/q(\alpha_n|x) \vee h\log\alpha_n \vee \varepsilon(q(\alpha_n|x)|x)) \to 0$ as $n\to \infty$,
\item $\beta_n/\alpha_n\to 0$ as $n\to \infty$,
\end{itemize}
then
$$
\frac{\sigma_n^{-1}(x)}{\log(\alpha_n/\beta_n)}\left(\frac{\hat q_n^{\mbox{\tiny W},\phi}(\beta_n|x)} {q(\beta_n|x)}-1\right)\tod {\mathcal{N}}(0,V_\phi(x)).
$$
\end{Coro}
As an example, in the case of the functional Hill and Pickands estimators, we obtain
\begin{eqnarray*}
V_{\phi_1}(x)&=& \gamma^2(x)\left(\sum_{j=1}^J \frac{2(J-j) +1}{\tau_j} - J^2\right) \left /
\left( \sum_{j=1}^J \log (1/\tau_j)\right)^2\right..\\
V_{\phi_{\mbox{\tiny FP}}}(x)&=& \frac{\gamma^2(x) (2^{2\gamma(x)+1}+1)}{4(\log 2)^2(2^{\gamma(x)}-1)^2}.
\end{eqnarray*}
Clearly, $V_{\phi_{\mbox{\tiny FP}}}(x)$ is
the variance of the classical Pickands estimator, see for instance~\cite{deHaanFer}, Theorem~3.3.5. 

\section{Illustration on simulated data}
\label{simul}

\begin{figure}
\begin{center}
\includegraphics[width=10cm,height=10cm]{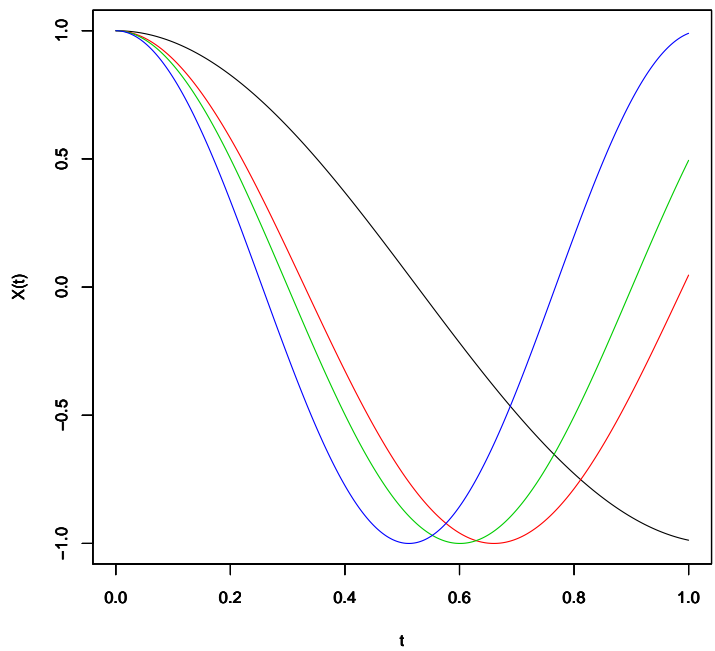} 
\end{center}
\caption{Four realizations of the random function $X(.)$.}
\label{realisations}
\end{figure}

The finite sample performance is illustrated on $N=200$ replications of a sample of size $n=500$
from a random pair $(X,Y)$, where the functional covariate $X \in E=L^2[0,1]$ is defined by
$X(t)=\cos (2\pi Z t)$ for all $t \in [0,1]$
where $Z$ is uniformly distributed on $[1/4,1]$. 
Some examples of simulated random functions $X$ are depicted on Figure~\ref{realisations}.
Besides, the conditional distribution of $Y$ given $X$ is a Burr distribution 
(see Table~1) with parameters $\tau(X)=2$ and 
$\lambda(X) = 2/(8 \|X\|_2^2-3)$ with 
$$\|X\|_2^2 = \int_0^1 X^2(t)dt= \frac{1}{2} \left ( 1+\frac{\sin (4 \pi Z)}{4 \pi Z}\right ).$$
We focus on the estimation of $q(\beta_n|x)$ with $\beta_n= 5/n$. 
To this end, the functional Weissman estimator $\qweis$ is used with a piecewise linear kernel
$K(t) = (1.9-1.8t) \I \{t \in [0,1] \}$ and the triangular kernel $Q'$.
The conditional tail index is estimated by the functional Hill estimator
${\hat{\gamma}}_n^{\phi_1}$.
The choice of the semi-metric $d$ is a recurrent issue in functional estimation~(see~\cite{LivreToulouse}, Chapter~3).
Here, two semi-metrics are considered. The first one is defined for all $(s,t)\in E^2$
by $d_X(s,t)=\|s-t\|_2$ and coincides with the $L_2$ distance between functions.
Remarking that the conditional quantile $q(\alpha_n|X)$ depends only on $\|X\|^2_2$,
or equivalently on $Z$, another interesting semi-metric is
$d_Z(s,t)=\left | \|s\|_2^2-\|t\|_2^2 \right |$. Finally, in Section~\ref{realdata}, an example
of the use of a metric based on second derivatives is presented.

With such choices, the functional Weissman estimator $\qweis$ depends on three
parameters $h$, $\lambda$ and $\alpha_n$ and on the $\tau_j$'s used to compute ${\hat{\gamma}}_n^{\phi_1}$.


\noindent - The smoothing parameter $h$ is selected using the cross-validation approach introduced in~\cite{Yao}
and implemented for instance in~\cite{Test2,Statmed}:
$$
h^{opt}=\arg\min \left\{\sum^n_{i=1}\sum^n_{j=1}\left(\I{\{Y_i\ge Y_j\}}-{\habaF}_{,-i}(Y_j|X_i)\right)^2,\ h\in {\cal H} \right\}
$$
where ${\habaF}_{,-i}$ is the estimator (depending on $h$) given in~(\ref{defestproba}) computed
from the sample $\{(X_\ell,Y_\ell),~1 \le \ell \le n,~\ell \ne i\}$.
Here, ${\cal H}$ is a regular grid, ${\cal H}=\{h_1\leq h_2\leq \dots \leq h_M\}$ with
$h_1=1/100$, $h_M=1/10$ and $M=20$. Let us note that this approach was originally
proposed for finite dimensional covariates. Up to our knowledge, its optimality (with respect to
the mean integrated squared error for instance) is not established in the functional framework.
We refer to~\cite{Ferraty3} for such a work in functional regression.
\\
- In our experiments, the choice of the bandwidth $\lambda$ appeared to be less crucial than 
the other smoothing parameter $h$. It could have been selected with the same criteria as 
previously, but for simplicity reasons, it has been fixed to $\lambda=0.1$. \\
- The choice of $\alpha_n$ is equivalent to the choice of the number of upper order
statistics in the non-conditional extreme-value theory. It is still an open question,
even though some techniques have been proposed, see for instance~\cite{Daniel} for
a bootstrap based method. \\
- The selection of the $\tau_j$'s is equivalent to the selection of an estimator for the conditional tail index. Once again, extreme-value theory does not provide optimal solution to this problem.

\begin{figure}[p]
\begin{center}
\includegraphics[width=9cm]{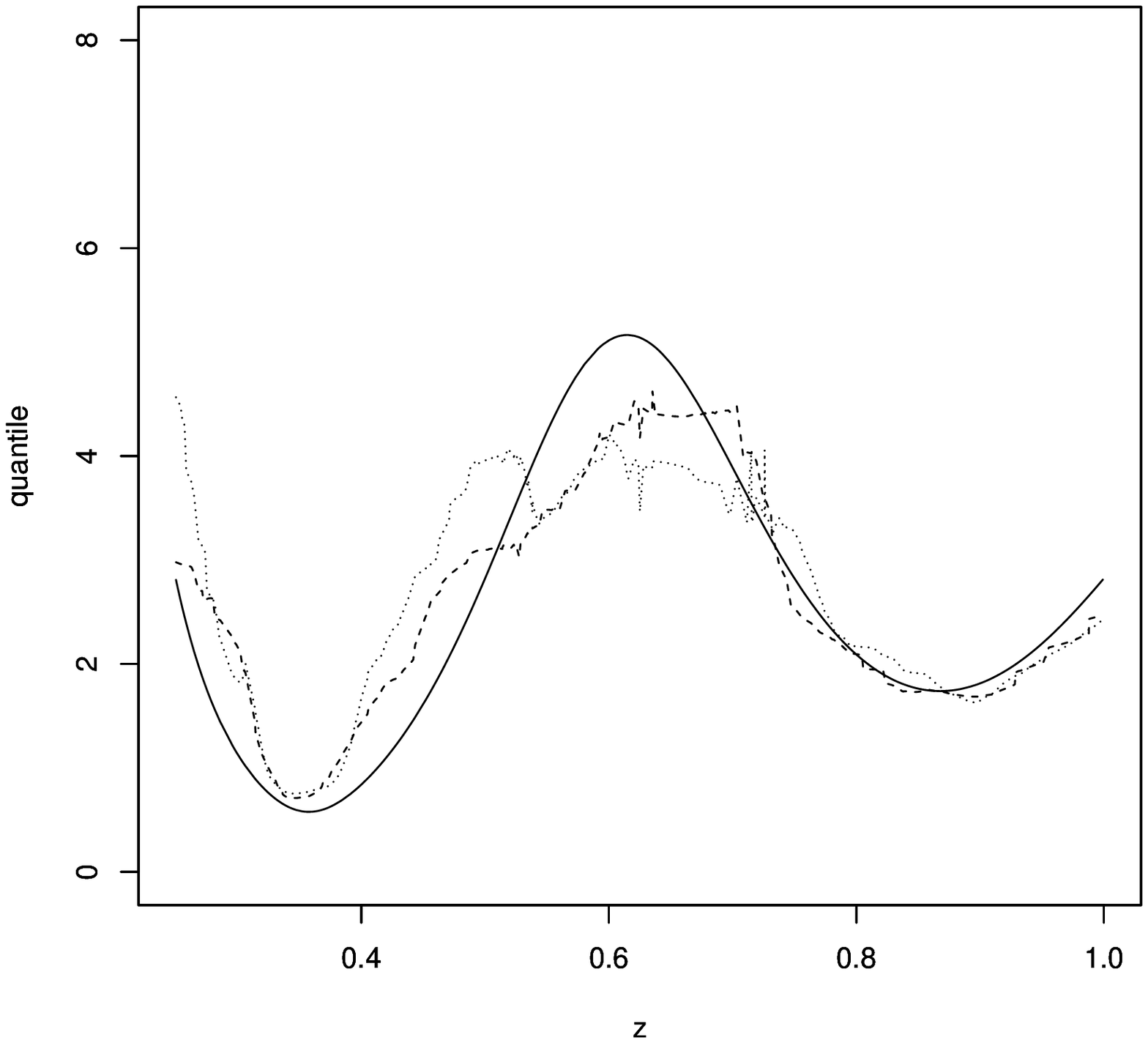} 
\includegraphics[width=9cm]{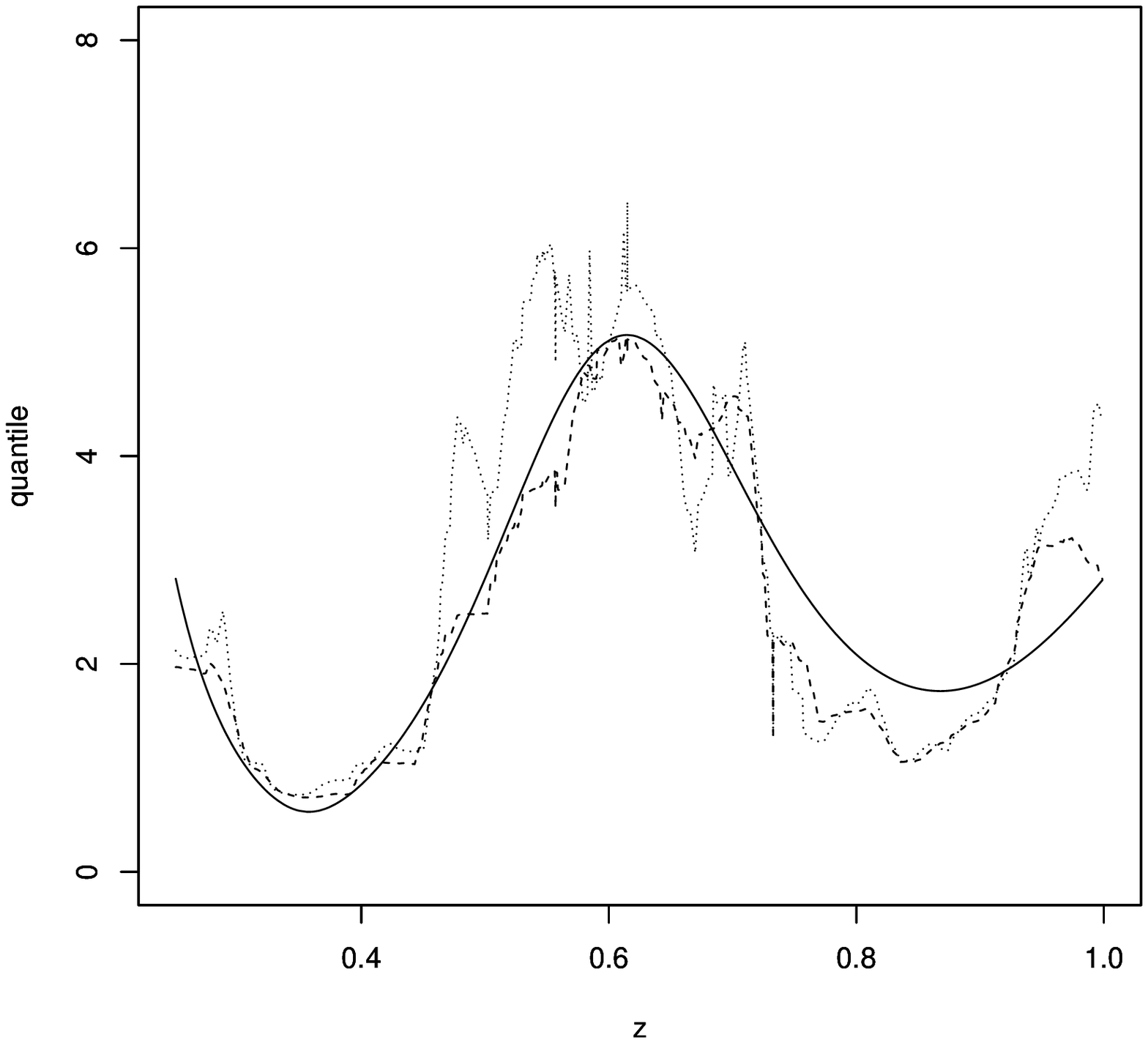} 
\end{center}
\caption{Comparison of the estimated quantile $\qweis$ corresponding to the 10\% quantile of the $L_2$-errors $\Delta_d^{(r)}$ with the true quantile function (continuous line).
Horizontally: $Z$, vertically: quantiles. Two sets of ($\alpha_n$,$\tau_j$) are considered: ($\alpha_n=15\log(n)/n$, $\tau_j=(1/j)^3$, dashed line) and 
($\alpha_n=10\log(n)/n$, $\tau_j=(1/j)^2$, dotted line). 
Top: semi-metric $d_Z$, bottom: semi-metric $d_X$. }
\label{quantiles10}
\end{figure}

\begin{figure}[p]
\begin{center}
\includegraphics[width=9cm]{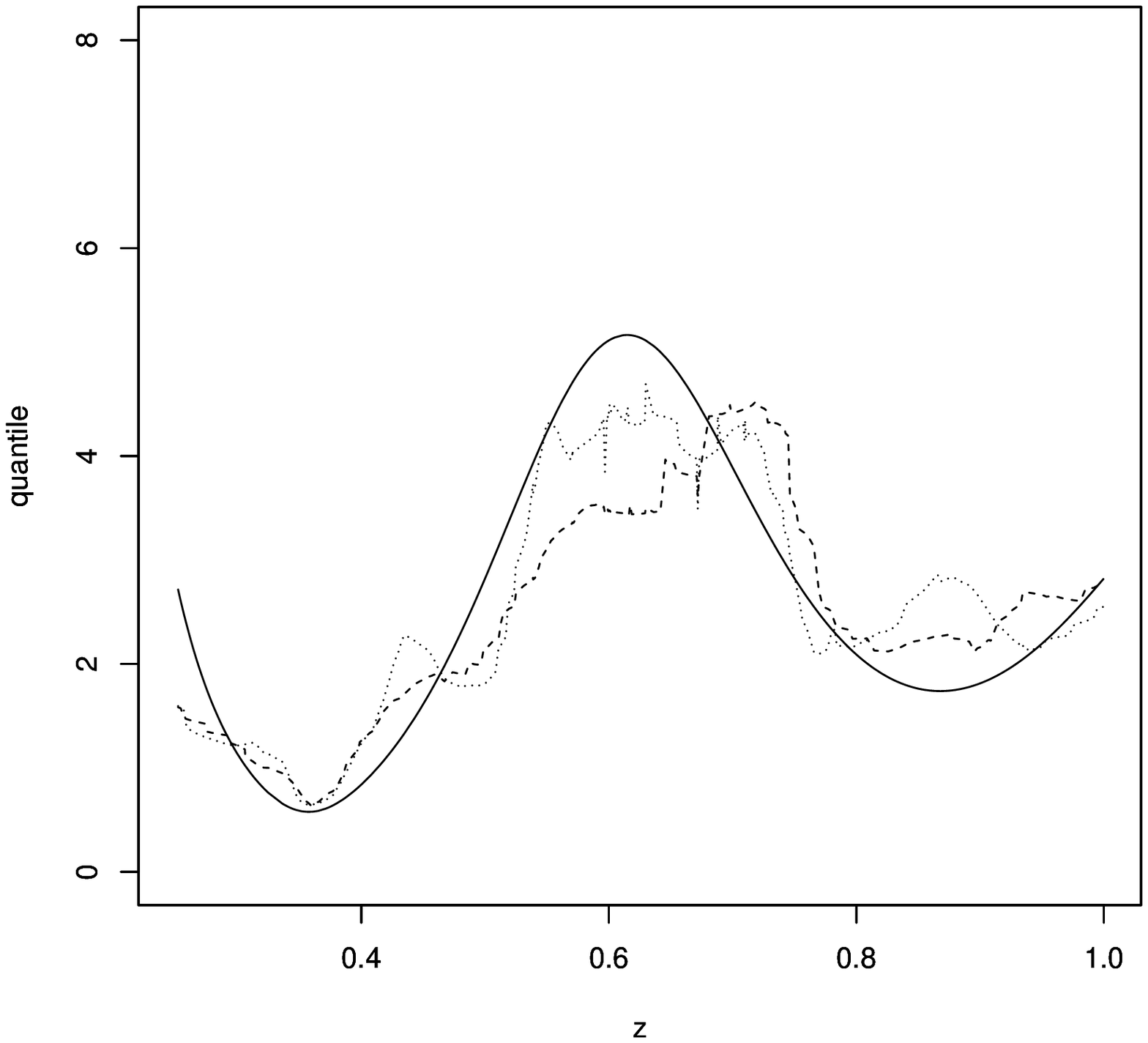} 
\includegraphics[width=9cm]{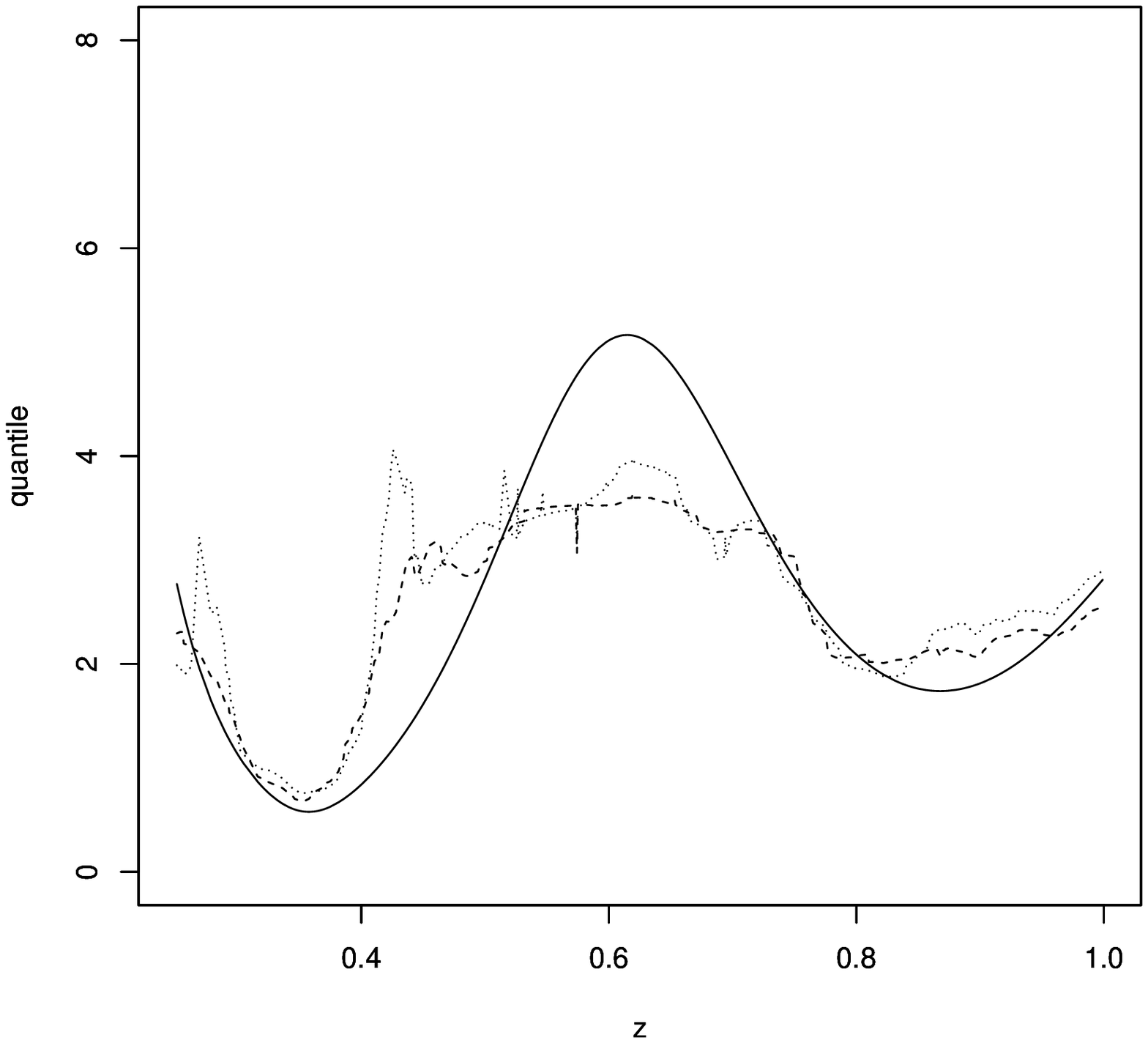} 
\end{center}
\caption{Comparison of the estimated quantile $\qweis$ corresponding to the median of the $L_2$-errors $\Delta_d^{(r)}$ with the true quantile function (continuous line).
Horizontally: $Z$, vertically: quantiles. Two sets of ($\alpha_n$,$\tau_j$) are considered: ($\alpha_n=15\log(n)/n$, $\tau_j=(1/j)^3$, dashed line) and 
($\alpha_n=10\log(n)/n$, $\tau_j=(1/j)^2$, dotted line). 
Top: semi-metric $d_Z$, bottom: semi-metric $d_X$. }
\label{quantiles50}
\end{figure}

\begin{figure}[p]
\begin{center}
\includegraphics[width=9cm]{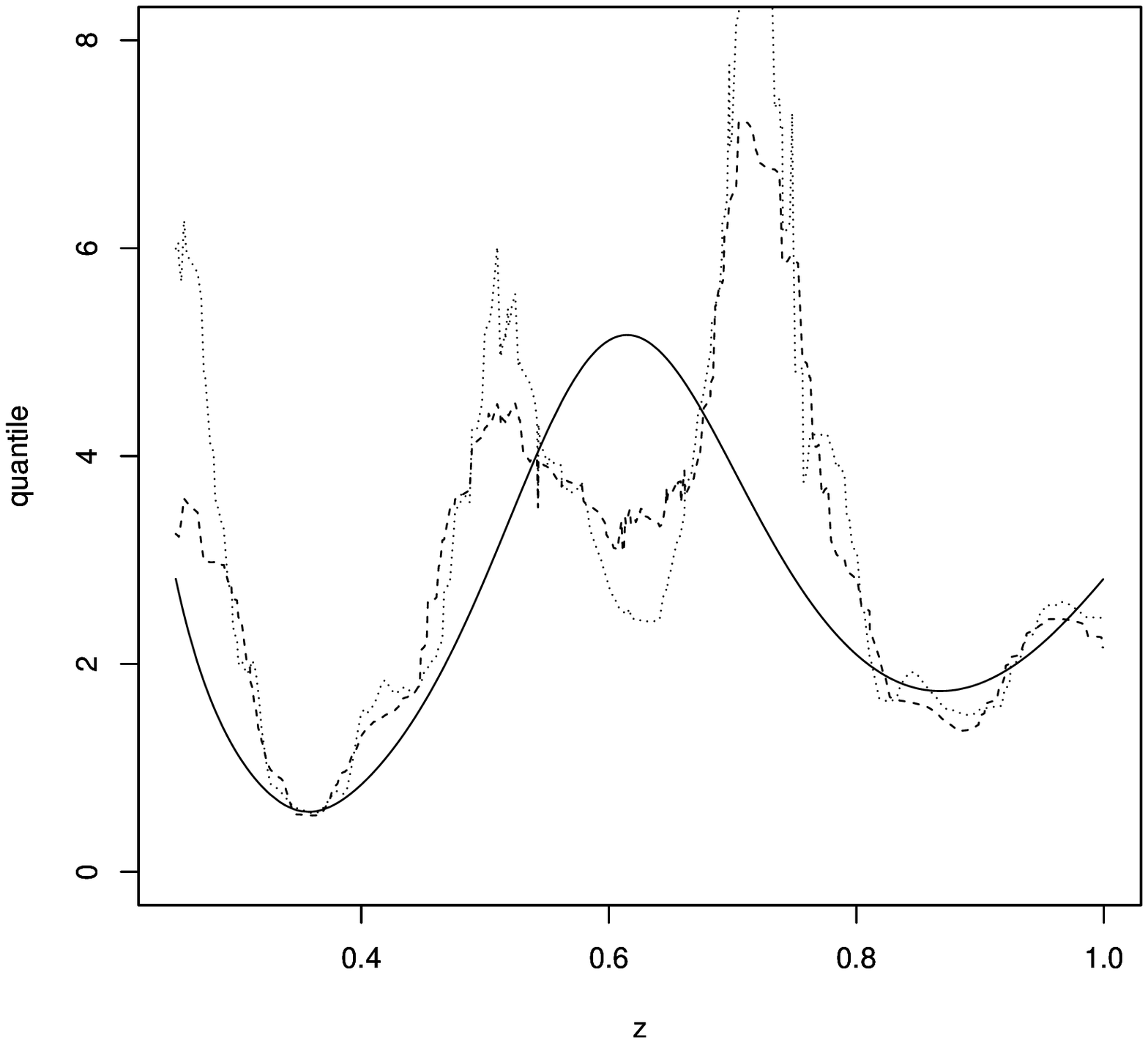} 
\includegraphics[width=9cm]{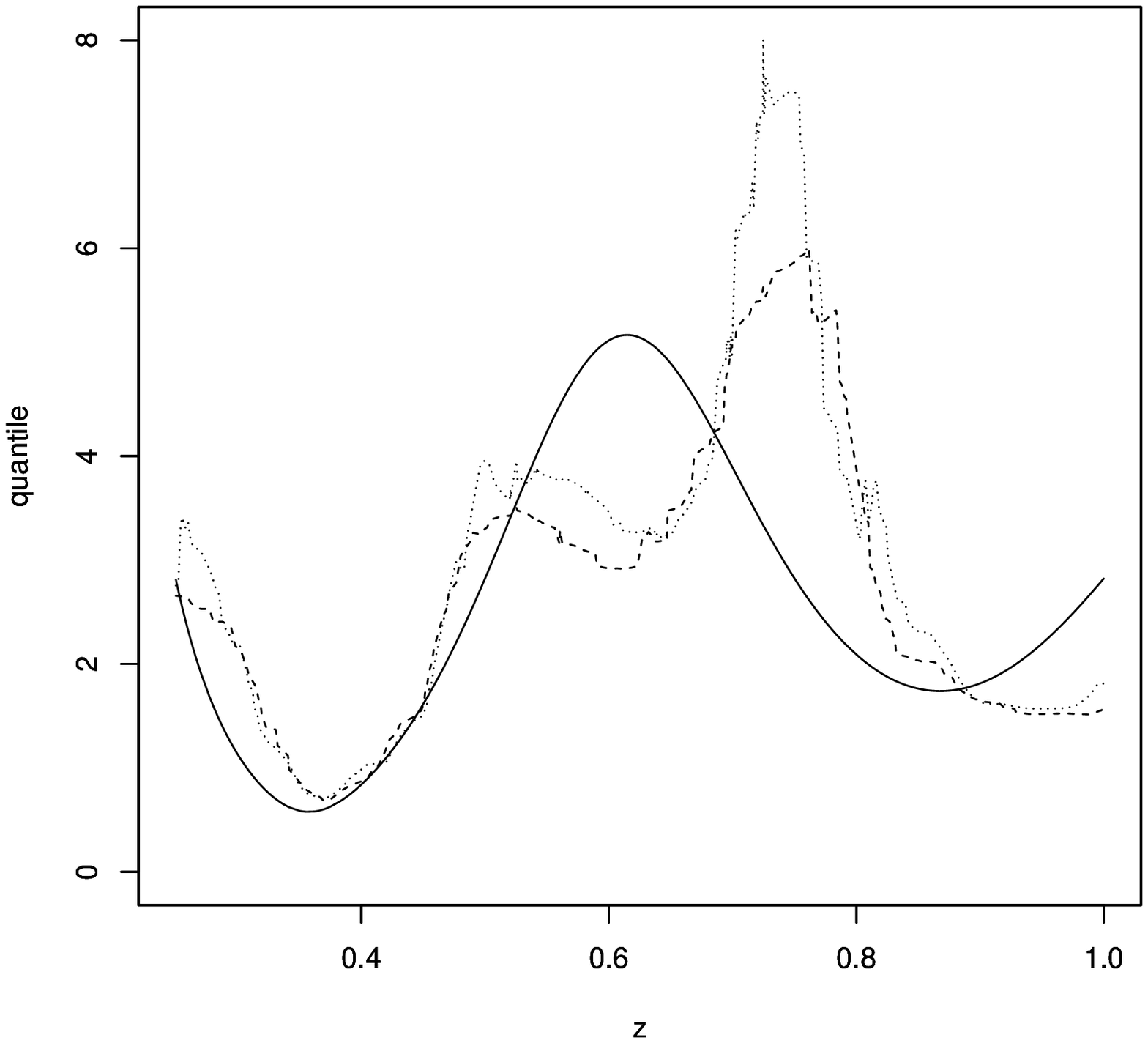} 
\end{center}
\caption{Comparison of the estimated quantile $\qweis$ corresponding to the 90\% quantile of the $L_2$-errors $\Delta_d^{(r)}$ with the true quantile function (continuous line).
Horizontally: $Z$, vertically: quantiles. Two sets of ($\alpha_n$,$\tau_j$) are considered: ($\alpha_n=15\log(n)/n$, $\tau_j=(1/j)^3$, dashed line) and 
($\alpha_n=10\log(n)/n$, $\tau_j=(1/j)^2$, dotted line). 
Top: semi-metric $d_Z$, bottom: semi-metric $d_X$. }
\label{quantiles90}
\end{figure}

In order to assess the impact of the choice of $\alpha_n$ and $\tau_j$'s, the $L_2$-errors 
$$
\Delta_d^{(r)} = \sum_{i=1}^n \left ( \hat q_n^{\mbox{\tiny W}}(\beta_n|X_i)^{(r)}-q(\beta_n|X_i) \right )^2,
$$
$r=1,\dots,N$ have been computed.
Here, $\hat q_n^{\mbox{\tiny W}}(\beta_n|X_i)^{(r)}$ is the estimation computed on the $r$th 
replication and $d$ can be either $d_X$ or $d_Z$.
Different values of $\alpha_n$ and $\tau_j$ are investigated: $\alpha_n=c\log (n)/n$ with $c \in \{5,10,15,20\}$ and $\tau_j=(1/j)^s$ with $s \in \{1,2,3,10\}$. The median, 10\% quantile and 90\%
quantile of the $\Delta_d^{(r)}$ errors 
are collected in Table~\ref{alphatau}. For a fixed value of $s$, the best error obtained with the semi-metric $d_Z$ is always smaller than the best error obtained with $d_X$ (both displayed in bold font).
Let us note that the optimal value of $c$ does not seem to depend on the semi-metric.
Besides, it will appear in the following that the estimations are not, at least visually, very sensitive with respect to the choice of $\alpha_n$ (or equivalently $c$) and $\tau_j$ (or equivalently $s$).
In Figure~\ref{quantiles10}--\ref{quantiles90},  the estimator $\qweis$ is represented as a function of $Z$. 
The estimator has been computed for two sets of ($\alpha_n$, $\tau_j$): 
($\alpha_n=15\log(n)/n$, $\tau_j=(1/j)^3$) and ($\alpha_n=10\log(n)/n$, $\tau_j=(1/j)^2$) and for the two semi-metrics $d_X$ and $d_Z$.
We limited ourselves to the representation of the estimator computed on the replications giving rise to the median, 10\% quantile and 90\% quantile of the $L_2$-errors $\Delta_d^{(r)}$, $r=1,\dots,N$. 
It appears that there is no
visual significative difference between the two choices of ($\alpha_n$, $\tau_j$). 

\section{Illustration on real data}
\label{realdata}

In this section, we propose to illustrate the behaviour of our large conditional quantiles estimators on functional chemometric data.  It concerns $n=215$ samples of finely chopped meat (see for example \cite{FV02} for more details). For each unit $i$ taken among this sample, we observe one spectrometric curve 
$\chi_i$ discretized at 100 wavelengths $\lambda_1,\ldots,\lambda_{100}$.
The covariate $x_i$ is thus defined by $x_i=(x_{i,1},\dots,x_{i,100})^t$
with $x_{i,j}=\chi_i(\lambda_j)$ for all $j=1,\dots,100$.
Each variable $x_{i,j}$ is the $-\log_{10}$ of the transmittance
recorded by the Tecator Infratec Food and Feed Analyzer spectrometer.
The dataset can be found at {\tt{http://lib.stat.cmu.edu/datasets/tecator}}.

Clearly, the covariate $x_i$ is in fact a discretized curve but, as mentioned in~\cite{LMS93}, the fineness of the grid spanning the discretization allows us to consider each subject as a continuous curve. Hence, the covariate can be considered as belonging to an infinite dimensional space $E$.
For each spectrometric curve $\chi_i$, the fat content ${\tilde{Y}}_i\in[0,100]$ (in percentage) is given. 
Since these values are bounded they cannot satisfy model~\Au~and we propose to use as variable of interest the inverse of the fat content defined as: $Y_i=100/{\tilde{Y}}_i\in[1,\infty)$, $i=1,\ldots,n$. 

In the following, the semi-metric based on the second derivative is 
adopted, as advised in~\cite{LivreToulouse}, Chapter~9:
$$
d^2(\chi_i,\chi_j)= \int \left ( \chi_i^{(2)}(t) - \chi_j^{(2)}(t) \right )^2dt,
$$
where $\chi^{(2)}$ denotes the second derivative of $\chi$. 
To compute this semi-metric, one can use an approximation of the functions $\chi_i$ and $\chi_j$ based on B-splines as proposed in~\cite{LivreToulouse}, Chapter~3. Here, we limit ourselves to a discretized version $\tilde d$ of $d$:
$$
\tilde d^2(x_i,x_j) = \sum_{l=2}^{99} \left\{ (x_{i,l+1} - x_{j,l+1}) + 
 (x_{i,l-1} - x_{j,l-1}) - 2 (x_{i,l}-x_{j,l}) \right\}^2.
$$
Other semi-metrics could be considered: Functional Principal Component Analysis (FPCA) or Multivariate
Partial Least-Squares Regression (MPLSR) are useful tools for computing proximities between curves
in reduced dimensional spaces, see~\cite{LivreToulouse}, Section~3.4.

\begin{figure}[p]
\begin{center}
\epsfig{figure=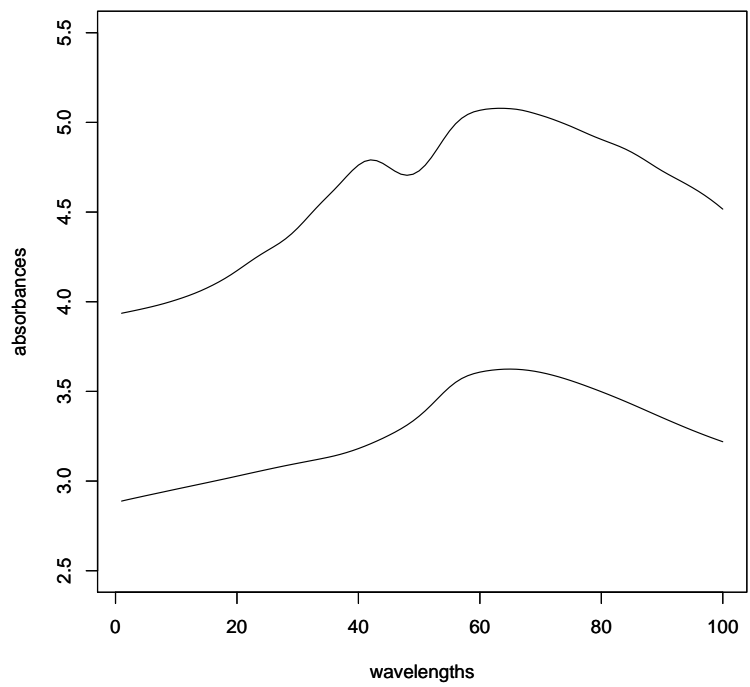,width=0.7\textwidth}
\caption{Selected spectrometric curves $\chi_{i_0}$ and $\chi_{i_1}$. }
\label{figcurves}
\end{center}
\end{figure}

We propose to estimate the large conditional quantile of order $\beta_n=5/n$ 
in a given direction of the space $E$. 
More precisely, we focus on the segment
$[\chi_{i_0},\chi_{i_1}]$ where $\chi_{i_0}$ and $\chi_{i_1}$ denote the most
different curves in the sample, {\it i.e.}
$$
(i_0,i_1)=\arg\max_{1\leq i<j\leq n}  \tilde d(x_i,x_j).
$$
The selected curves $\chi_{i_0}$ and $\chi_{i_1}$ are plotted in
Figure~\ref{figcurves}. Since these curves appear to be smooth, the chosen semi-metric, which
is based on the second derivative, seems to be well adapted.
The conditional quantile to estimate is
$q(\beta_n, t(\xi))$ where $t(\xi)= \xi\chi_{i_1} + (1-\xi) \chi_{i_0}$ for $\xi\in[0,1]$.
To this end, the functional Weissman estimator is considered with the same kernels as in
the previous section. The selected smoothing parameters are $h=0.02$ and $\lambda=0.1$.

\begin{figure}[p]
\begin{center}
\includegraphics[width=8cm]{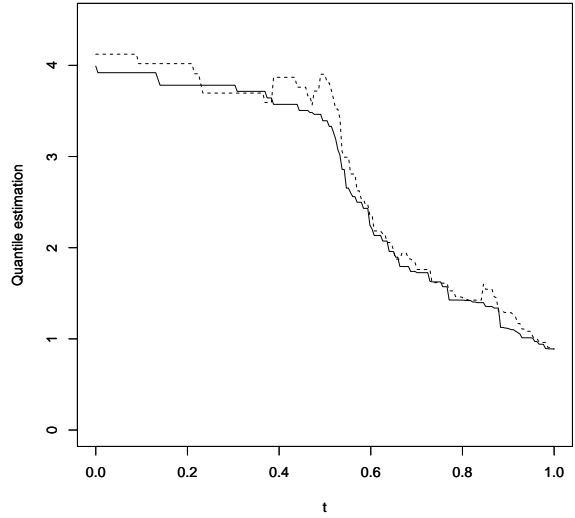} 
\end{center}
\caption{Quantile estimate of order $\beta_n=5/n$ as a function of $t(\xi)=\xi\chi_{i_1}+(1-\xi)\chi_{i_0}$,
 $\xi\in[0,1]$. Continuous line: $\tau_j=(1/j)^3$ and $\alpha_n=15\log(n)/n$, dashed line: $\tau_j=(1/j)^2$ and $\alpha_n=10\log(n)/n$.}
\label{quant1et2}
\end{figure}

\begin{figure}[p]
\begin{center}
\includegraphics[width=8cm]{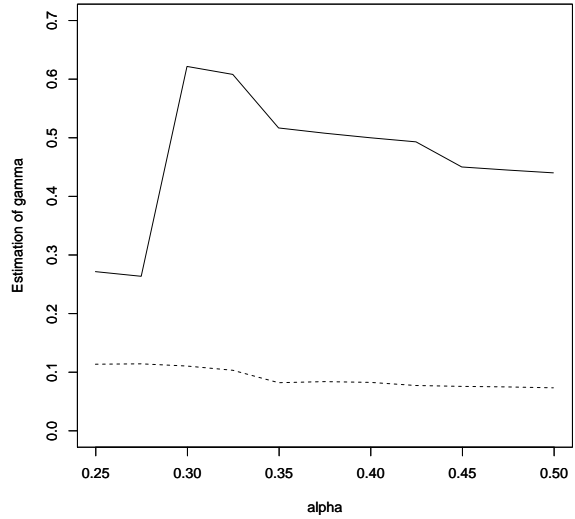} 
\end{center}
\caption{Estimation of the conditional tail-index $\gamhill$ as a function of $\alpha_n$ with $\tau_j=(1/j)^2$.
Continuous line: $x=\chi_{i_0}$, dashed line: $x=\chi_{i_1}$.}
\label{gam1et2}
\end{figure}

The estimated quantile $\hat q_n^{\mbox{\tiny W}}(\beta_n,t(\xi))$
is plotted as a function of $\xi$ in Figure~\ref{quant1et2} for
different values of weights $\tau_j$ and probability $\alpha_n$.
Here again, it appears that the estimated quantiles are not too sensitive
with respect to these parameters.
The globally decreasing shape
of the curves indicates that heaviest tails 
({\it i.e.} largest values of $\gamma(t(\xi))$)
are found in the neighbourhood of the curve $\chi_{i_0}$
({\it i.e.} for small values of~$\xi$). At the opposite, lightest tails
are found in the neighbourhood of the curve $\chi_{i_1}$.
These results are confirmed by Figure~\ref{gam1et2}: The estimated conditional
tail-index $\gamhill$ is larger for $x=\chi_{i_0}$ than for $x=\chi_{i_1}$.
These very different shapes confirm a strong heterogeneity of the sample in
terms of tail behaviour.

\section{Further work}
\label{further}

Our further work will consist in establishing uniform convergence results.
The rate of uniform strong consistency of the csf estimator $\habaF(y|x)$
defined in~(\ref{defestproba})
is already known since~\cite{Ferraty2} for fixed $y$.
The first step will then to adapt this result for $y=y_n\to\infty$ as $n\to\infty$.
On this basis, it should be possible to 
get uniform results for $\hat q(\alpha_n|x)$ (see~(\ref{defestquant}))  in the case of large conditional quantiles
withing the sample, {\it ie.}  $\alpha_n\to 0$ with $n \varphi_x(x) \alpha_n\to \infty$.
The last step would be to extend these results to $\qweis$ defined in~(\ref{defweis})
when $\beta_n\to0$ arbitrarily fast. Such results would require the uniform convergence
of $\hat\gamma_n(x)$, the estimator of the conditional tail index.

\section{Appendix: Proofs}
\label{proofs}

\subsection{Preliminary results}

\noindent The following two lemmas are of analytical nature. 
The first one is dedicated to the control of the local variations
of the csf 
when the quantity of interest $y$ goes to infinity.

\begin{Lem}
\label{lemlip}
Let $x\in E$ and suppose \Au~and \Ad~hold. \\
(i) If $y_n\to\infty$ and $h\log y_n\to 0$
as $n\to \infty$, then, for $n$ large enough,
$$
\sup_{x'\in B(x,h)} \left| \frac{\baF(y_n|x)}{\baF(y_n|x')}-1 \right|
\leq 2(\kappa_c + \kappa_\gamma+\kappa_\varepsilon)h\log y_n.
$$
(ii) If $y_n\to\infty$ and $y_n'\to \infty$
as $n\to \infty$, then, for $n$ large enough,
$$
\sup_{x'\in B(x,h)} \left| \frac{\baF(y'_n|x')}{\baF(y_n|x')}-1 \right|
\leq \left| \left(\frac{y_n}{y'_n}\right)^{2/\gamma(x)} -1\right|.
$$
\end{Lem}
\proof (i) Assumption \Au~yields, for all $x'\in B(x,h)$:
\begin{eqnarray*}
\left|\log\left( \frac{\baF(y_n|x)}{\baF(y_n|x')} \right)\right|
&\leq& \left|\log c(x) - \log c(x')\right| + \int_1^{y_n} \left( \left| \frac{1}{\gamma(x)}-\frac{1}{\gamma(x')} \right| + \left| \varepsilon(u|x)-\varepsilon(u|x') \right|
\right) \frac{du}{u} \\
&\leq& \kappa_c h +\int_1^{y_n} (\kappa_\gamma + \kappa_\varepsilon) h \frac{du}{u} \\
&\leq & (\kappa_c + \kappa_\gamma+\kappa_\varepsilon) h\log y_n ,
\end{eqnarray*}
eventually, from \Ad. Thus,
$$
\sup_{d(x,x')\leq h} \left|\log\left( \frac{\baF(y_n|x)}{\baF(y_n|x')} \right)\right|
=O(h\log y_n) \to 0
$$
as $n\to\infty$ and taking account of $\log(u+1)\sim u$ as $u\to 0$
gives the result.\\

\noindent (ii) Let us assume for instance $y'_n>y_n$. From \Au~we have
\begin{equation}
 \label{majo}
 \left| \frac{\baF(y'_n|x')}{\baF(y_n|x')}-1 \right|=
1 - \left(\frac{y'_n}{y_n}\right)^{-1/\gamma(x')} \exp\left( \int_{y_n}^{y'_n} 
\frac{\varepsilon(u|x')}{u}du \right)
\leq 1 - \left(\frac{y'_n}{y_n}\right)^{-1/\gamma(x')-|\varepsilon(y_n|x')|}.
\end{equation}
Now, $x'\in B(x,h)$ and \Ad~imply for $n$ large enough that 
$$
\frac{1}{\gamma(x')} + |\varepsilon(y_n|x')| \leq \frac{1}{\gamma(x)} + (\kappa_\varepsilon+\kappa_\gamma) h + |\varepsilon(y_n|x)| \leq \frac{2}{\gamma(x)}.
$$
Replacing in~(\ref{majo}), it follows that
$$
 \left| \frac{\baF(y'_n|x')}{\baF(y_n|x')}-1 \right| \leq 1 - \left(\frac{y'_n}{y_n}\right)^{-2/\gamma(x)}.
$$
The case $y'_n\leq y_n$ is similar.
\CQFD
\noindent The second lemma provides
a second order asymptotic expansion of the quantile function.
It is proved in~\cite{Test2}.
\begin{Lem}
\label{lemDLq}
Suppose \Au~hold.\\
(i) Let $0<\beta_n<\alpha_n$ with $\alpha_n\to0$ as $n\to\infty$. Then,
$$
|\log q(\beta_n|x) - \log q(\alpha_n|x) + \gamma(x)\log(\beta_n/\alpha_n)|=O(\log(\alpha_n/\beta_n)\varepsilon(q(\alpha_n|x)|x)).
$$
(ii) If, moreover, $\lim\inf \beta_n/\alpha_n>0$, then
$$
\frac{\beta_n^{\gamma(x)}q( \beta_n|x)}{\alpha_n^{\gamma(x)}q(\alpha_n|x)} 
=1+ O(\varepsilon(q(\alpha_n|x)|x)).
$$
\end{Lem}
 
\noindent The following lemma provides a control on the moments $\mu_x^{(\tau)}(h)$ for all $\tau>0$,
the case $\tau=1$ being studied in~\cite{LivreToulouse}, Lemma~4.3. The proof is straightforward.
\begin{Lem}
\label{lemmu}
Suppose \At~holds. For all $\tau>0$ and $x\in E$,
$
0<C_1^\tau \varphi_x(h) \leq \mu_x^{(\tau)}(h)\leq C_2^\tau \varphi_x(h).
$
\end{Lem}
\noindent The following lemma provides a geometrical interpretation
of the condition $n \varphi_x(h) \baF(y_n|x)\to\infty$.
\begin{Lem}
\label{leminter}
Suppose \Au, \Ad~hold and let $y_n\to\infty$ such that
$h\log y_n\to 0$ as $n\to \infty$. Consider the subset of $E\times \R$ defined as
$R_n(x)=B(x,h)\times (y_n,\infty)$ where $x\in E$ is such that $\varphi_x(h)> 0$.
Then,
 $\PP(\exists i\in\{1,\dots,n\}, (X_i,Y_i)\in R_n(x))\to 1$
as $n\to\infty$ if, and only if, $n \varphi_x(h) \baF(y_n|x)\to\infty$.
\end{Lem}
\proof Since $(X_i,Y_i)$, $i=1,\dots,n$ are independent and identically distributed random variables,
\begin{equation}
\label{eqproba}
 \PP(\exists i\in\{1,\dots,n\}, (X_i,Y_i)\in R_n(x))=
1- (1 - \PP((X,Y)\in R_n(x)))^n
\end{equation}
where
\begin{eqnarray*}
\PP((X,Y)\in R_n(x)))&=& \E(\I\{X \in B(x,h) \cap Y\geq y_n \}) \\
&=& \E(\I\{X \in B(x,h)  \} \baF(y_n|X))\\
&= & \baF(y_n|x)\varphi_x(h) +  \baF(y_n|x)\E\left( \left(\frac{ \baF(y_n|X) }{ \baF(y_n|x)} -1\right) \I\{X \in B(x,h)  \}\right).
\end{eqnarray*}
In view of Lemma~\ref{lemlip}(i), we have
$$
\E\left( \left|\frac{ \baF(y_n|X) }{ \baF(y_n|x)} -1\right| \I\{X \in B(x,h)  \}\right)\leq 2(\kappa_c + \kappa_\gamma + \kappa_\varepsilon) \varphi_x(h) h \log y_n 
$$
and therefore
$$
\PP((X,Y)\in R_n(x))= \baF(y_n|x)\varphi_x(h) (1 + O(h\log y_n)).
$$
Clearly, this probability converges to 0 as $n\to\infty$ and 
thus~(\ref{eqproba}) can be rewritten as
$$
\PP(\exists i\in\{1,\dots,n\}, (X_i,Y_i)\in R_n(x))=1- 
\exp\left( -  n \varphi_x(h) \baF(y_n|x)(1+o(1))\right),
$$
which converges to 1 if and only if $ n \varphi_x(h) \baF(y_n|x)\to\infty$.
\CQFD

\noindent Let us remark that the kernel estimator~(\ref{defestproba})
can be rewritten as $\habaF(y|x)=\hat\psi_n(y,x)/\hat g_n(x)$
with
\begin{eqnarray*}
\hat \psi_n(y,x)&=& \frac{1}{n \muu} \sum_{i=1}^n K(d(x,X_i)/h) Q((Y_i-y)/\lambda),\\
\hat g_n(x) &=&  \frac{1}{n \muu} \sum_{i=1}^n K(d(x,X_i)/h).
\end{eqnarray*}
Lemma~\ref{lemdensite} and  Lemma~\ref{lempsi} are respectively dedicated to the asymptotic properties of $\hat g_n(x)$ and $\hat\psi_n(y,x)$.
\begin{Lem}
\label{lemdensite}
Suppose \At~holds and let $x\in E$ such that $\varphi_x(h)>0$. We have:
\begin{description}
\item[(i)] $\E(\hat g_n(x)) = 1$.
\item[(ii)] If, moreover, $\varphi_x(h)\to 0$ as $h\to 0$ then
$$
0<\liminf n \varphi_x(h) \operatorname{var}(\hat g_n(x)) \leq \limsup n \varphi_x(h) 
\operatorname{var}(\hat g_n(x))<\infty.
$$
\end{description}
\end{Lem}
Therefore, under \At, if $\varphi_x(h)\to 0$ and $n\varphi_x(h)\to\infty$ then 
$\hat g_n(x)$ converges to $1$ in probability.

\noindent \proof {\bf (i)} is straightforward.\\
{\bf (ii)} Standard calculations yields
$$
n \varphi_x(h) \mbox{var}(\hat g_n(x))=  \varphi_x(h)\left( \frac{\mu_x^{(2)}(h)}{(\mu_x^{(1)}(h))^2}-1\right)
$$
and Lemma \ref{lemmu} entails
$$
(C_1/C_2)^2\leq \varphi_x(h) \frac{\mu_x^{(2)}(h)}{(\mu_x^{(1)}(h))^2} \leq (C_2/C_1)^2.
$$
The condition $\varphi_x(h)\to 0$ concludes the proof.\CQFD
\newpage
\begin{Lem}
\label{lempsi}
Suppose \Au~-- \Aq~hold. Let $x\in E$ such that $\varphi_x(h)>0$ and introduce $y_{n,j}=a_j y_n(1+o(1))$ for $j=1,\dots, J$ with
$0<a_1<a_2< \dots < a_J$ and where $J$ is a positive integer.
If $y_n\to\infty$ such that $h\log y_n\to 0$, $\lambda/y_n\to 0$ and $n\varphi_x(h) \baF(y_n|x) \to \infty$ as $n\to \infty$,
then
\begin{description}
\item [(i)] $\E(\hat\psi_n(y_{n,j},x)) = \baF(y_{n,j}|x) (1+ O(h \log y_n \vee \lambda/y_n))$,
for $j=1,\dots, J$.
\item [(ii)] The random vector 
$$
\left\{ \Lambda_n^{-1}(x) \left(\frac{\hat\psi_n(y_{n,j},x)
-\E(\hat\psi_n(y_{n,j},x))}{\baF(y_{n,j}|x)}\right)\right\}_{j=1,\dots,J}
$$
is asymptotically Gaussian, centered, with covariance matrix 
$C(x)$ where $C_{j,j'}(x) = a^{1/\gamma(x)}_{j\wedge j'}$ 
for $(j,j')\in \{1,\dots, J\}^2$.
\end{description}
\end{Lem}
\proof {\bf (i)} The $(X_i,Y_i)$, $i=1,\dots,n$ being identically
distributed, we have
\begin{eqnarray*}
\E(\hat\psi_n(y_{n,j},x))&=&\frac{1}{\muu} \E\{K(d(x,X)/h)Q((Y-y_{n,j})/\lambda)\}\\
&=& \frac{1}{\muu} \E\{K(d(x,X)/h) \E(Q((Y-y_{n,j})/\lambda)|X)\}
\end{eqnarray*}
Taking account of \Aq, it follows that
$$
\E(Q((Y-y_{n,j})/\lambda)|X)= \baF(y_{n,j}|X)
+ \int_{-1}^1 Q'(u) (\baF(y_{n,j}+\lambda u|X)-\baF(y_{n,j}|X))du
$$
and thus the bias can be expanded as
\begin{equation}
\label{decbias}
\E(\hat\psi_n(y_{n,j},x))- \baF(y_{n,j}|x)=: T_{1,n} + T_{2,n},
\end{equation}
where we have defined
\begin{eqnarray*}
T_{1,n}&=& \frac{1}{\muu} \E\{K(d(x,X)/h)(\baF(y_{n,j}|X)-\baF(y_{n,j}|x))\},\\
T_{2,n}&=&\frac{1}{\muu} \E\left\{K(d(x,X)/h)\baF(y_{n,j}|X) \int_{-1}^1 Q'(u) \left(\frac{\baF(y_{n,j}+\lambda u|X)}{\baF(y_{n,j}|X)}-1\right)du \right\}.
\end{eqnarray*}
Focusing on $T_{1,n}$ and taking account of \At, it follows that
\begin{eqnarray*}
T_{1,n}
&=& \frac{1}{\muu} \E(K(d(x,X)/h)(\baF(y_{n,j}|X)-\baF(y_{n,j}|x))\I\{d(x,X)\leq h\})\\
&=&\frac{\baF(y_{n,j}|x)}{\muu} \E\left(K(d(x,X)/h)\left(\frac{\baF(y_{n,j}|X)}{\baF(y_{n,j}|x)}-1\right)\I\{d(x,X)\leq h\}\right).
\end{eqnarray*}
Lemma~\ref{lemlip}(i) implies that
$$
 \left| \frac{\baF(y_{n,j}|X)}{\baF(y_{n,j}|x)}-1 \right|\I\{d(x,X)\leq h\}
\leq 2(\kappa_c + \kappa_\gamma + \kappa_\varepsilon)h\log y_{n,j} \leq 3 (\kappa_c + \kappa_\gamma + \kappa_\varepsilon)h \log y_n,
$$
eventually and therefore
\begin{equation}
\label{eqT1}
|T_{1,n}| = \baF(y_{n,j}|x) O(h\log y_n).
\end{equation}
Let us now consider $T_{2,n}$. From Lemma~\ref{lemlip}(ii), for all $u\in[-1,1]$, we eventually have
$$
\left|\frac{\baF(y_{n,j}+\lambda u|X)}{\baF(y_{n,j}|X)}-1\right|\I\{d(x,X)\leq h\}
\leq \left| \left( 1+ \frac{\lambda u}{y_{n,j}}\right)^{2/\gamma(x)} -1\right|
 \leq C_{\gamma(x)} \frac{\lambda}{y_{n,j}},
$$
since $\lambda/y_n\to 0$ as $n\to\infty$ and where $C_{\gamma(x)}$ is a positive constant. As a consequence,
\begin{eqnarray}
\nonumber|T_{2,n}|& \leq &C_{\gamma(x)} \frac{\lambda}{y_{n,j}} \frac{1}{\muu} \E(K(d(x,X)/h)\baF(y_{n,j}|X))\\
& = &C_{\gamma(x)}\frac{\lambda}{y_{n,j}} 
(\baF(y_{n,j}|x)+ T_{1,n})= \baF(y_{n,j}|x)O(\lambda/y_n)
\label{eqT2}
\end{eqnarray}
in view of (\ref{eqT1}). Collecting (\ref{decbias}), (\ref{eqT1}) and (\ref{eqT2})
concludes the first part of the proof.\\
\noindent {\bf (ii)} Let $\beta\neq 0$ in $\R^J$
and consider the random variable
$$
\Psi_n=\sum_{j=1}^J \beta_j \left(\frac{\hat\psi_n(y_{n,j},x)
-\E(\hat\psi_n(y_{n,j},x))}{\Lambda_n(x)\baF(y_{n,j}|x)}\right) 
=: \sum_{i=1}^n Z_{i,n},
$$
where, for all $i=1,\dots,n$, the random variable $Z_{i,n}$ is defined by
\begin{eqnarray*}
n \Lambda_n(x)\muu Z_{i,n}&=&\left\{  \sum_{j=1}^J \frac{\beta_j K(d(x,X_i)/h) Q((Y_i-y_{n,j})/\lambda) }{\baF(y_{n,j}|x)}\right. \\
&- &\left.\E\left( \sum_{j=1}^J \frac{\beta_j K(d(x,X_i)/h)Q((Y_i-y_{n,j})/\lambda) }{\baF(y_{n,j}|x)} \right)\right\}.
\end{eqnarray*}
Clearly, $\{Z_{i,n},\; i=1,\dots,n\}$ is a set of centered, independent 
and identically distributed random variables. Let us determine an asymptotic expansion of their variance:
\begin{eqnarray}
\nonumber
\mbox{var}(Z_{i,n}) &=& \frac{1}{n^2 (\muu)^2 \Lambda_n^2(x)} \mbox{var}
\left( \sum_{j=1}^J \beta_j K(d(x,X_i)/h) 
\frac{Q((Y_i-y_{n,j})/\lambda) }{\baF(y_{n,j}|x)} \right)\\
\nonumber
&=& \frac{1}{n^2 (\muu)^2\Lambda_n^2(x)} \beta^t B(x) \beta\\
\label{varzin}
&=& \frac{\baF(y_n|x)}{n \mud} \beta^t B(x) \beta,
\end{eqnarray}
where $B(x)$ is the $J\times J$ covariance matrix 
with coefficients defined for $(j,j')\in\{1,\dots,J\}^2$ by
\begin{eqnarray*}
B_{j,j'}(x)&=& \frac{A_{j,j'}(x)}{\baF(y_{n,j}|x)\baF(y_{n,j'}|x)},\\
A_{j,j'}(x) &=& \mbox{cov}\left\{ K(d(x,X)/h) Q((Y-y_{n,j})/\lambda)
, \; K(d(x,X)/h) Q((Y-y_{n,j'})/\lambda)\right\}\\
&=& \E\left\{ K^2(d(x,X)/h)Q((Y-y_{n,j})/\lambda) Q((Y-y_{n,j'})/\lambda)\right\}\\
&-& \E\{ K(d(x,X)/h) Q((Y-y_{n,j})/\lambda) \} \E\{  K(d(x,X)/h)Q((Y-y_{n,j'})/\lambda)  \}\\
&=:& T_{3,n} - T_{4,n}.
\end{eqnarray*}
Let us first focus on $T_{3,n}$:
\begin{equation}
\label{eqT3}
 T_{3,n}= \E\{K^2(d(x,X)/h) \E(Q((Y-y_{n,j})/\lambda)Q((Y-y_{n,j'})/\lambda)|X)\}
\end{equation}
and remark that
$$
\E(Q((Y-y_{n,j})/\lambda)Q((Y-y'_{n,j})/\lambda)|X)=:\Omega(y_{n,j},y_{n,j'}) + \Omega(y_{n,j'},y_{n,j})
$$
where we have defined
\begin{eqnarray*}
\Omega(y,z)&=&\frac{1}{\lambda} \int_\R Q'((t-y)/\lambda)Q((t-z)/\lambda) \baF(t|X) dt\\
&=& \int_{-1}^1 Q'(u) Q(u + (y-z)/\lambda) \baF(y+u\lambda|X)du.
\end{eqnarray*}
Let us consider the case $j< j'$. We thus have $a_j< a_{j'}$ and consequently
$(y_{n,j}-y_{n,j'})/\lambda\to -\infty$ as $n\to\infty$. Therefore, for $n$ large enough
$u+ (y_{n,j}-y_{n,j'})/\lambda<-1 $ and $Q(u+ (y_{n,j}-y_{n,j'})/\lambda)=0$.
It follows that, eventually 
$\Omega(y_{n,j},y_{n,j'})=0$. Similarly, for $n$ large enough
 $Q(u+ (y_{n,j'}-y_{n,j})/\lambda)=1$ and 
$$
\Omega(y_{n,j'},y_{n,j})= \int_{-1}^1 Q'(u) \baF(y_{n,j'}+u\lambda|X)du.
$$
For symmetry reasons, it follows that, for all $j\neq j'$,
$$
\E(Q((Y-y_{n,j})/\lambda)Q((Y-y'_{n,j})/\lambda)|X)=\int_{-1}^1 Q'(u) \baF(y_{n,j\vee j'}+u\lambda|X)du= \E(Q((Y-y_{n,j\vee j'})/\lambda)|X),
$$
and replacing in (\ref{eqT3}) yields
$$
 T_{3,n}= \E\{K^2(d(x,X)/h)\E(Q((Y-y_{n,j\vee j'})/\lambda)|X)\}=
\E\{K^2(d(x,X)/h)Q((Y-y_{n,j\vee j'})/\lambda)\}.
$$
Now, since $K^2$ is a kernel also satisfying assumption \At,
part {\bf (i)} of the proof implies
\begin{equation}
\label{result}
T_{3,n} = \mud \baF(y_{n,j\vee j'}|x)(1+ O(h \log y_n \vee \lambda/y_n)),
\end{equation}
for all $j\neq j'$. 
In the case where $j=j'$, by definition,
$$
T_{3,n}= \E\{K^2(d(x,X)/h)\E(Q^2((Y-y_{n,j})/\lambda)|X)\}
$$
where $K^2$ is a kernel also satisfying assumption \At~and where the
pdf associated to $Q^2$ satisfies assumption \Aq. Consequently,
(\ref{result}) also holds for $j=j'$.
Second, part {\bf (i)} of the proof implies
$$
T_{4,n} = (\muu)^2 \baF(y_{n,j}|x)\baF(y_{n,j'}|x) (1+ O(h \log y_n\vee\lambda/y_n)).
$$
As a consequence, 
\begin{eqnarray*}
A_{j,j'}(x) &=& \mud \baF(y_{n,j\vee j'}|x) (1+O(h\log y_n \vee \lambda/y_n))\\
&-& (\muu)^2 \baF(y_{n,j}|x)\baF(y_{n,j'}|x) (1+O(h\log y_n \vee \lambda/y_n))
\end{eqnarray*}
leading to 
$$
B_{j,j'}(x) = \frac{\mud}{\baF(y_{n,j\wedge j'}|x)} \left(1+O(h\log y_n \vee\lambda/y_n)
- \frac{(\muu)^2}{\mud} \baF(y_{n,j\wedge j'}|x)  (1+O(h\log y_n \vee\lambda/y_n))\right).
$$
In view of Lemma~\ref{lemmu}, $(\muu)^2/\mud$ is bounded and taking account of
$\baF(y_{n,j\wedge j'}|x)\to 0$ as $n\to\infty$ yields
$$
B_{j,j'}(x)=  \frac{\mud}{\baF(y_{n,j\wedge j'}|x)}(1+o(1)).
$$
Now, from the regular variation property~(\ref{defl}), it is easily seen that
$$
\baF(y_{n,j\wedge j'}|x) = a_{j\wedge j'}^{-1/\gamma(x)}\baF(y_n|x)(1+o(1))
$$
entailing
$
B_{j,j'}(x) =  C_{j,j'}(x)\mud/\baF(y_n|x) (1+o(1)).
$
Replacing in~(\ref{varzin}), it follows that
$$
\mbox{var}(Z_{i,n}) = \frac{\beta^t C(x) \beta}{n} (1+o(1)),
$$
for all $i=1,\dots,n$. As a preliminary conclusion, 
var$(\Psi_n)\to \beta^t C(x) \beta$ as $n\to\infty$.
Consequently, Lyapounov criteria for the asymptotic normality of sums of
triangular arrays reduces to
$
\sum_{i=1}^{n}\E\left| Z_{i,n}\right|^{3}
=n\E\left| Z_{1,n}\right| ^{3}\rightarrow0
$
as $n\to\infty$.
Next, remark that $Z_{1,n}$ is a bounded random variable:
\begin{eqnarray*}
|Z_{1,n}| &\leq &\frac{ 2 C_2 \sum_{j=1}^J |\beta_j|}{n \Lambda_n(x) \muu \baF(y_{n,J}|x)} \\
&=& 2 C_2 a_J^{1/\gamma(x)} \frac{\muu}{\mud}\sum_{j=1}^J |\beta_j|\Lambda_n(x)(1+o(1))\\
&\leq &2 (C_2/C_1)^2 a_J^{1/\gamma(x)}\sum_{j=1}^J |\beta_j|\Lambda_n(x)(1+o(1));
\end{eqnarray*}
in view of Lemma~\ref{lemmu} and thus, 
\begin{eqnarray*}
n\E\left| Z_{1,n}\right| ^{3} &\leq& 2 (C_2/C_1)^2 a_J^{1/\gamma(x)} \sum_{j=1}^J |\beta_j|  \Lambda_n(x)  n{\mbox{var}}(Z_{1,n})(1+o(1))\\
& =&  2 (C_2/C_1)^2  a_J^{1/\gamma(x)}\sum_{j=1}^J |\beta_j| \beta^t C(x) \beta \Lambda_n(x)(1+o(1))\to 0
\end{eqnarray*}
as $n\to\infty$ in view of Lemma~\ref{lemmu}. As a conclusion, $\Psi_n$ converges in distribution
to a centered Gaussian random variable with variance $\beta^t C(x) \beta$ for all $\beta\neq 0$ in $\R^J$. The result is proved.\CQFD

\subsection{Proofs of main results}

\paragraph{Proof of Theorem \ref{thproba}.}

Keeping in mind the notations of Lemma~\ref{lempsi}, the following expansion holds
\begin{equation}
\label{term123}
\Lambda_n^{-1}(x)\sum_{j=1}^J \beta_j \left(\frac{\habaF(y_{n,j}|x)}{\baF(y_{n,j}|x)} -1\right)
=:\frac{\Delta_{1,n}+\Delta_{2,n}-\Delta_{3,n}}{\hat g_n(x)},
\end{equation}
where 
\begin{eqnarray*}
\Delta_{1,n}&=&\Lambda_n^{-1}(x) \sum_{j=1}^J \beta_j \left(\frac{\hat \psi_n(y_{n,j},x) -
\E(\hat \psi_n(y_{n,j},x)) }{\baF(y_{n,j}|x)} \right)\\
\Delta_{2,n}&=& \Lambda_n^{-1}(x)\sum_{j=1}^J \beta_j \left(\frac{ \E(\hat\psi_n(y_{n,j},x)) -  \baF(y_{n,j}|x) }{\baF(y_{n,j}|x)} \right)\\
\Delta_{3,n}& =& \left( \sum_{j=1}^J \beta_j\right)\Lambda_n^{-1}(x)\left( \hat g_n(x)-1 \right).
\end{eqnarray*}
Let us highlight that assumptions $nh^2\varphi_x(h) \log^2(y_n) \baF(y_n|x) \to 0$ and $n \varphi_x(h)\baF(y_n|x) \to \infty$ imply that $h\log y_n\to 0$ as $n\to \infty$.
Thus, from Lemma~\ref{lempsi}(ii), the random term $\Delta_{1,n}$ can be
rewritten as 
\begin{equation}
\label{term1}
\Delta_{1,n}=\sqrt{\beta^t C(x) \beta} \xi_n,
\end{equation}
where $\xi_n$ converges to a standard Gaussian random variable.
The nonrandom term $\Delta_{2,n}$ is controlled with Lemma~\ref{lempsi}(i):
\begin{equation}
\label{term2}
\Delta_{2,n}= O(\Lambda_n^{-1}(x)( h\log y_n \vee \lambda/y_n))=o(1).
\end{equation}
Finally, $\Delta_{3,n}$ can be bounded by Lemma~\ref{lemdensite} and Lemma~\ref{lemmu}:
\begin{equation}
\label{term3}
\Delta_{3,n}=  O_P(\Lambda_n^{-1}(x)(n\varphi_x(h))^{-1/2})=
 O_P(\baF(y_n|x))^{1/2}=o_P(1).
\end{equation}
Collecting (\ref{term123})--(\ref{term3}), it follows that
$$
\hat g_n(x) \Lambda_n^{-1}(x)\sum_{j=1}^J \beta_j \left(\frac{\habaF(y_{n,j}|x)}{\baF(y_{n,j}|x)} -1\right)
= \sqrt{\beta^t C(x) \beta} \xi_n +o_P(1).
$$
Finally, $\hat g_n(x)\toP 1$ concludes the proof. \CQFD

\paragraph{Proof of Theorem \ref{thquant}.}

Introduce for $j=1,\dots,J$, 
\begin{eqnarray*}
\alpha_{n,j}&=&\tau_j \alpha_n,\\
\sigma_{n,j}(x)& =&q(\alpha_{n,j}|x) \sigma_n(x),\\
v_{n,j}(x)&=&\alpha_{n,j}^{-1}\gamma(x) \sigma_n^{-1}(x),\\
W_{n,j}(x)&=&v_{n,j}(x)\left(\habaF(q(\alpha_{n,j}|x) + \sigma_{n,j}(x) z_j|x)-\baF(q(\alpha_{n,j}|x) + \sigma_{n,j}(x) z_j|x)\right),\\
a_{n,j}(x)&=&v_{n,j}(x)\left(\alpha_{n,j}-\baF(q(\alpha_{n,j}|x) + \sigma_{n,j}(x) z_j |x)\right),
\end{eqnarray*}
and  $z_j\in\R$. Let us study the asymptotic behavior of $J$-variate function
defined by
$$
\Phi_n(z_1,\dots,z_J)=\PP\left(\bigcap_{j=1}^J \left\{
\sigma_{n,j}^{-1}(x) (\hat q_n(\alpha_{n,j}|x)-q(\alpha_{n,j}|x))\leq z_j
\right\}\right)
= \PP\left(\bigcap_{j=1}^J \left\{
W_{n,j}(x) \leq a_{n,j}(x) \right\}\right).
$$
We first focus on the nonrandom term $a_{n,j}(x)$.
Under \Au, $\baF(.|x)$ is differentiable.
Thus, for all $j\in\{1,\dots,J\}$ there exists $\theta_{n,j}\in(0,1)$
such that 
\begin{equation}
\label{DL}
\baF(q(\alpha_{n,j}|x)|x) - 
\baF(q(\alpha_{n,j}|x) + \sigma_{n,j}(x) z_j|x)=
-  \sigma_{n,j}(x) z_j \baF'(q_{n,j}|x),
\end{equation}
where $q_{n,j}=q(\alpha_{n,j}|x) + \theta_{n,j}\sigma_{n,j}(x) z_j$. 
It is clear that $q(\alpha_{n,j}|x)\to\infty$ and $\sigma_{n,j}(x)/q(\alpha_{n,j}|x)\to0$
as $n\to\infty$. As a consequence, $q_{n,j}\to\infty$ and thus
\Au~entails
\begin{equation}
\label{RV}
\lim_{n\to\infty} \frac{q_{n,j} \baF'(q_{n,j}|x)}{ \baF(q_{n,j}|x)} = -1/\gamma(x).
\end{equation}
Moreover, since $q_{n,j}= q(\alpha_{n,j}|x)(1+o(1))$ 
and $\baF(.|x)$ is regularly varying at infinity, it follows that
$\baF(q_{n,j}|x)=\baF( q(\alpha_{n,j}|x)|x)(1+o(1))=\alpha_{n,j}(1+o(1))$.
In view of~(\ref{DL}) and~(\ref{RV}), we end up with
\begin{equation}
\label{aj}
a_{n,j}(x)=
  \frac{ v_{n,j}(x)\sigma_{n,j}(x) \alpha_{n,j} z_j}{\gamma(x)q(\alpha_{n,j}|x)}(1+o(1))
=  z_j (1+o(1)).
\end{equation}
Let us now turn to the random term $W_{n,j}(x)$. 
Defining $a_j=\tau_j^{-\gamma(x)}$, $y_{n,j}=q(\alpha_{n,j}|x) + \sigma_{n,j}(x) z_j$ for $j=1,\dots,J$ and $y_n=q(\alpha_n|x)$,
we have $y_{n,j}=q(\alpha_{n,j}|x)(1+o(1))= a_j y_n(1+o(1))$ since $q(.|x)$ is
regularly varying at 0 with index $-\gamma(x)$.
Using the same argument, it is easily shown that
$\log y_n=  -\gamma(x)\log(\alpha_n)(1+o(1))$. 
As a consequence, Theorem~\ref{thproba} applies and the random vector
$$
\left\{ \frac{\sigma_n^{-1}(x)}{v_{n,j}(x) \baF(y_{n,j}|x)} W_{n,j} \right\}_{j=1,\dots,J} = (1+o(1)) \left\{  \frac{W_{n,j}}{\gamma(x)} \right\}_{j=1,\dots,J}
$$
converges to a centered Gaussian random variable with covariance
matrix $C(x)$.
Taking account of~(\ref{aj}), we obtain that $\Phi_n(z_1,\dots,z_J)$ 
converges to the cumulative distribution function of a centered Gaussian distribution with covariance
matrix $\gamma^2(x)C(x)$ evaluated at $(z_1,\dots,z_J)$,
which is the desired result.
\CQFD

\paragraph{Proof of Theorem \ref{thweis}.} The proof is based on the following expansion:
$$
\frac{\sigma_n^{-1}(x)}{\log(\alpha_n/\beta_n)}(\log(\qweis)-\log(q(\beta_n|x))) = \frac{\sigma_n^{-1}(x)}{\log(\alpha_n/\beta_n)} (Q_{n,1}+Q_{n,2}+Q_{n,3})
$$
where we have introduced
\begin{eqnarray*}
 Q_{n,1}&=&  \sigma_n^{-1}(x)(\hat\gamma_n(x)-\gamma(x)),\\
 Q_{n,2}&=& \frac{\sigma_n^{-1}(x)}{\log(\alpha_n/\beta_n)} \log (\hat q_n(\alpha_n|x)/q(\alpha_n|x)),\\
 Q_{n,3}&=& \frac{\sigma_n^{-1}(x)}{\log(\alpha_n/\beta_n)} (\log q(\alpha_n|x) - \log q(\beta_n|x) + \gamma(x)\log(\alpha_n/\beta_n)).
\end{eqnarray*}
First, $Q_{n,1}\tod {\mathcal{N}}(0,V(x))$ as a straightforward consequence of the assumptions.
Second, Theorem~\ref{thquant} implies that $\hat q_n(\alpha_n|x)/q(\alpha_n|x)\toP 1$ and 
$$
Q_{n,2}=\frac{\sigma_n^{-1}(x)}{\log(\alpha_n/\beta_n)} \left(\frac{\hat q_n(\alpha_n|x)}{q(\alpha_n|x)}-1\right)(1+o_P(1))=\frac{O_P(1)}{\log(\alpha_n/\beta_n)}.
$$
Consequently, $Q_{n,2}\toP 0$ as $n\to\infty$. Finally, from Lemma~\ref{lemDLq}(i),
$
Q_{n,3}= O( \sigma_n^{-1}(x)\varepsilon(q(\alpha_n|x)|x)),
$
which converges to 0 in view of the assumptions.
\CQFD

\paragraph{Proof of Theorem \ref{thgamma}.} The following expansion holds for all $j=1,\dots,J$:
\begin{equation}
\label{eqdev1}
\log \hat q_n(\tau_j \alpha_n|x) = \log q(\alpha_n|x) 
+ \log\left(\frac{ q(\tau_j\alpha_n|x) }{ q(\alpha_n|x)}\right)
+ \log\left(\frac{ \hat q_n(\tau_j\alpha_n|x) }{ q(\tau_j\alpha_n|x)}\right).
\end{equation}
First, Lemma~\ref{lemDLq}(ii) entails that
\begin{equation}
\label{eqdev2}
\log\left(\frac{ q(\tau_j\alpha_n|x) }{ q(\alpha_n|x)}\right)=\gamma(x)\log(1/\tau_j)+ O(\varepsilon(q(\alpha_n|x)|x)),
\end{equation}
where the $O(\varepsilon(q(\alpha_n|x)|x))$ is not necessarily uniform in $j=1,\dots,J$.
Second, it follows from Theorem~\ref{thquant} that
\begin{equation}
\label{eqdev3}
\log\left(\frac{ \hat q_n(\tau_j\alpha_n|x) }{ q(\tau_j\alpha_n|x)}\right)= \sigma_n(x) \xi_{n,j}
\end{equation}
where $(\xi_{n,1},\dots,\xi_{n,J})^t$
converges to a centered Gaussian random vector with covariance
matrix $\gamma^2(x)\Sigma$.
Replacing~(\ref{eqdev2}) and~(\ref{eqdev3}) in~(\ref{eqdev1}) yields
$$
\log \hat q_n(\tau_j \alpha_n|x) = \log q(\alpha_n|x) + \gamma(x)\log(1/\tau_j) + \sigma_n(x) \xi_{n,j}+ O(\varepsilon(q(\alpha_n|x)|x)),
$$
for all $j=1,\dots,J$ and therefore, in view of the shift invariance property of $\phi$, we have
$$
\phi\left( \{\log \hat q_n(\tau_j \alpha_n|x)\}_{j=1,\dots,J}\right)= 
\phi\left( \{\gamma(x)\log(1/\tau_j) + \sigma_n(x) \xi_{n,j}+ O(\varepsilon(q(\alpha_n|x)|x)) \}_{j=1,\dots,J} \right).
$$
A first order Taylor expansion yields:
\begin{eqnarray*}
\phi\left( \{\log \hat q_n(\tau_j \alpha_n|x)\}_{j=1,\dots,J}\right)&= &
\phi\left( \gamma(x)v \right) + 
\sum_{j=1}^J (\sigma_n(x) \xi_{n,j}+ O(\varepsilon(q(\alpha_n|x)|x))) \frac{\partial \phi}{\partial x_j}(\gamma(x) v)\\
&+& O_P\left( \sum_{j=1}^J (\sigma_n(x) \xi_{n,j}+ O(\varepsilon(q(\alpha_n|x)|x)))^2\right).
\end{eqnarray*}
Thus, under the condition $\sigma^{-1}_n(x) \varepsilon(q(\alpha_n|x)|x)\to 0$ as $n\to\infty$, it follows that
$$
\sigma^{-1}_n(x) (\phi\left( \{\log \hat q_n(\tau_j \alpha_n|x)\}_{j=1,\dots,J}\right)-\phi\left( \gamma(x)v \right))
= \sum_{j=1}^J \xi_{n,j}  \frac{\partial \phi}{\partial x_j}(\gamma(x) v) + o_P(1).
$$
Taking into account of the scale invariance property of $\phi$, we finally obtain
$$
\sigma^{-1}_n(x) (\hat \gamma^{\phi}_n(x)-\gamma(x))
= \frac{1}{\phi( v)} \sum_{j=1}^J \xi_{n,j}  \frac{\partial \phi}{\partial x_j}(\gamma(x) v) + o_P(1)
$$
and the conclusion follows.
\CQFD



\newpage
\bibliographystyle{plain}

\begin{table}[hb]
\begin{center}
$
 \begin{array}{|l|c|c|c|c|}
\hline
 &\baF(y|x) & \gamma(x) & c(x) & \varepsilon(y|x) \\
\hline
&&&&\\
\mbox{Pareto} & y^{-\theta(x)} & \displaystyle\frac{1}{\theta(x)} & 1 & 0 \\

\mbox{Cauchy} & \displaystyle\frac{1}{\pi}\tan^{-1}(1/y) + \frac{1}{2}(1-\mbox{sign}(y)) & 1 &  \displaystyle\frac{1}{4} & \displaystyle\frac{2}{3} \frac{1}{y^2}(1+o(1)) \\

\mbox{Fr\'echet} & 1- \exp( -y^{-\theta(x)}) & \displaystyle\frac{1}{\theta(x)} &  1-e^{-1} & \displaystyle\frac{\theta(x)}{2}{y^{-\theta(x)}}(1+o(1)) \\

\mbox{Burr} & (1 +  y^{\tau(x)})^{-\lambda(x)} & \displaystyle\frac{1}{\lambda(x)\tau(x)} &  2^{-\lambda(x)} & \lambda(x)\tau(x){y^{-\tau(x)}}(1+o(1)) \\
&&&&\\
\hline
\end{array}
$
\end{center}
\caption{Examples of distributions satisfying~\Au. Their parameters $\theta(x)$, $\tau(x)$ and $\lambda(x)$ are positive.}
\label{exemples}
\end{table}

\begin{table}[hb]
\begin{center}
$
\begin{array}{|c c||c|c|c|c|}
\hline
 \ & \ & c=5 & c=10 & c=15 & c=20 \\
\hline
s=1 & d=d_X & 10^5 \; [2787,10^8] & 863 \; [363,3103] & {\bf 793} \; [311,2402] & 936 \; [344,2492] \\
 \ & d=d_Z  & 10^{7} \; [7208,10^{12}] & 860 \; [323,3137] & {\bf 688} \; [287,2242] & 751 \; [352,2586] \\
\hline
s=2 & d=d_X & 6391   \; [871,10^5]  & 429 \; [176,1347] & 349 \; [144,1056] & {\bf 341} \; [151,1106] \\
 \ & d=d_Z  & 10^{5} \; [2310,10^8] & 449 \; [195,1525] & 342 \; [156,1212] & {\bf 329} \; [148,1260] \\
\hline
s=3 & d=d_X & 2300  \; [570,10^5]  & 318 \; [126,1083] & 272 \; [111,792] & {\bf 231} \; [099,650] \\
 \ & d=d_Z  & 13651 \; [1436,10^6]  & 309 \; [141,1301] & 277 \; [115,863] & {\bf 228} \; [109,672]\\
\hline
s=10 & d=d_X & 430 \; [191,1963] & {\bf 392} \; [191,6357] & 665 \; [423,943] & 895 \; [633,1164]\\
 \ & d=d_Z   & 795 \; [328,8062] & {\bf 378} \; [170,6477] & 660 \; [372,933] & 894 \; [577,1155]\\
\hline
\end{array}
$
\end{center}
\caption{Median [10\% quantile, 90\% quantile] of the $L_2$-errors $\Delta_d^{(r)}$ for $\alpha_n =c\log (n)/n$, $c\in \{5,10,15,20\}$ and $\tau_j=(1/j)^s$, $s \in \{1,2,3,10\}$ for the two semi-metrics $d_X$ and $d_Z$.}
\label{alphatau}
\end{table}


\end{document}